\documentclass{amsart}
\usepackage{graphicx}
\usepackage{amssymb}
\usepackage{epstopdf}
\usepackage{nicefrac}
\usepackage{enumerate}
\usepackage{float}
\usepackage{tabularx}
\usepackage{hyperref}

\input xy
\xyoption {all}

\newtheorem{thm}{THEOREM}[section]

\newtheorem{cor}[thm]{COROLLARY}
\newtheorem{defn}[thm]{DEFINITION}

\newtheorem{lemma}[thm]{LEMMA}

\newtheorem{prop}[thm]{PROPOSITION}

 \newcommand{\ds}{\displaystyle}

\newcommand{\F}{{\mathcal F}}

\newcommand{\ab}{\alpha\beta}
\newcommand{\ba}{\beta\alpha}
\newcommand{\bB}{{\bf B}}
\newcommand{\bG}{{\bf G}}
\newcommand{\bH}{{\bf H}}
\newcommand{\cA}{{\mathcal A}}
\newcommand{\cB}{{\mathcal B}}
\newcommand{\cD}{{\mathcal D}}
\newcommand{\cE}{{\mathcal E}}
\newcommand{\cG}{{\mathcal G}}
\newcommand{\cI}{{\mathcal I}}
\newcommand{\cJ}{{\mathcal J}}
\newcommand{\cK}{{\mathcal K}}
\newcommand{\cP}{{\mathcal P}}
\newcommand{\cQ}{{\mathcal Q}}
\newcommand{\cR}{{\mathcal R}}

\newcommand{\cU}{{\mathcal U}}
\newcommand{\cV}{{\mathcal V}}
\newcommand{\cX}{{\mathcal X}}

\newcommand{\DF}{D_{\F}}
\newcommand{\Dv}{D^{\vec{v}}}
\newcommand{\dT}{{\bf d}_{\cX}}
\newcommand{\dX}{{\bf d}_{\fX}}
\newcommand{\E}{{\rm E}^+(\F)}
\newcommand{\SF}{{\rm S}(\F)}
\newcommand{\Ea}{{\rm E}^+_a(\F)}
\newcommand{\G}{\Gamma}
\newcommand{\g}{g_{\F}}
\newcommand{\GF}{{\cG}_{\F}}
\newcommand{\h}{h({\F})}

\newcommand{\hg}{{\bf g}}
\newcommand{\hh}{{\bf h}}

\newcommand{\mR}{{\mathbb R}}
\newcommand{\mS}{{\mathbb S}}
\newcommand{\mZ}{{\mathbb Z}}
\newcommand{\tcD}{{\widetilde{\cD}}}
\newcommand{\tcP}{{\widetilde{\cP}}}

\newcommand{\tcX}{{\widetilde{\cX}}}
\newcommand{\tU}{{\widetilde{U}}}
\newcommand{\tX}{{\widetilde{\fX}}}
\newcommand{\wh}{\widetilde{\bf h}} 

\newcommand{\oB}{{\overline{B}}}

\newcommand{\lF}{{\lambda_{\mathcal F}}} 
\newcommand{\lM}{{\lambda_{M}}} 
\newcommand{\e}{{\epsilon}} 
\newcommand{\eU}{{\epsilon_{\cU}}} 
\newcommand{\dF}{d_{\F}} 
\newcommand{\dFU}{\delta^{\F}_{\cU}} 
\newcommand{\whh}{{\widehat{\bf h}}}

\newcommand{\wtt}{{\widetilde t}}

\newcommand{\wtpi}{{\widetilde \pi}}
\newcommand{\wtvp}{{\widetilde \varphi}}
\newcommand{\vp}{{\varphi}}
\newcommand{\wtcP}{{\widetilde \cP}}

\newcommand{\fX}{{\mathfrak{X}}}

\begin{document}
\bibliographystyle{plain}

\title{Dynamics and the Godbillon-Vey Class of $C^1$ Foliations}

\author{Steven Hurder}
\address{Steven Hurder, Department of Mathematics, University of Illinois at Chicago, 322 SEO (m/c 249), 851 S. Morgan Street, Chicago, IL 60607-7045, USA}
\email{hurder@uic.edu}
\author{R\'emi Langevin}
\address{R\'emi Langevin, Universit\'e de Bourgogne - Franche Comt\'e, Institut de Math\'ematiques de Bourgogne UMR CNRS 5584, 21078 Dijon, France}
\email{Remi.Langevin@u-bourgogne.fr}
 \thanks{submitted March 27, 2004;  revised December 23, 2015}
 \date{}

\begin{abstract}
Let $\F$ be a codimension--one, $C^2$-foliation  on a manifold $M$ without boundary. In this work we show that if the Godbillon--Vey class $GV(\F) \in H^3(M)$ is non-zero, then $\F$ has a hyperbolic resilient leaf. Our approach is based on methods of $C^1$-dynamical systems, and does not use the classification theory of $C^2$-foliations.
 We first   prove that for a codimension--one $C^1$-foliation with non-trivial Godbillon measure, the set of infinitesimally expanding points $E(\F)$  has positive Lebesgue measure. We then prove that if $E(\F)$ has positive measure for a $C^1$-foliation, then $\F$ must have a hyperbolic resilient leaf, and hence its geometric entropy must be positive. The proof of this uses a pseudogroup version of the Pliss Lemma. The first statement then follows, as a $C^2$-foliation with non-zero Godbillon-Vey class has non-trivial    Godbillon measure. These results apply for both the case when $M$ is compact, and when $M$ is an open manifold.
\end{abstract}

\maketitle

 \vspace{-.2in}

 
   \tableofcontents

\vfill
\eject

\section{Introduction}\label{sec-intro}

 Godbillon and Vey  introduced in \cite{GodbillonVey1971}  the   invariant 
$GV(\F) \in H^3(M; \mR)$ named after them, which is defined for a codimension-one $C^2$-foliation $\F$ of a   manifold $M$ without boundary.
While the definition of the Godbillon-Vey class is elementary, 
  understanding    its relations to the geometric and dynamical properties of the foliation $\F$ remains an open problem. 
In the paper  \cite{Thurston1972} by Thurston, where he showed that the Godbillon-Vey class can assume   a continuous range of  values for foliations of closed $3$-manifolds, he also included   Figure~1, which illustrated the concept of ``hellical wobble'', which he suggested   gives a relation between the value of this class and the Riemannian geometry of the foliation. This geometric relation was made precise in a work by Reinhart and Wood \cite{ReinhartWood1973}.  More recently, Langevin and Walczak   in \cite{LW2008,Walczak2009,Walczak2013} gave  further    insights  into   the geometric meaning of the Godbillon-Vey invariant for smooth foliations of closed 3-manifolds, in terms of  the conformal geometry of the leaves of the foliation.

The Godbillon-Vey class   appears    in a surprising variety of contexts, such as the Connes-Moscovici work on the cyclic cohomology of Hopf algebras  \cite{Connes1986,ConnesMoscovici1998,Connes1994} which interprets the class in non-commutative geometry setting.  The   works of Leichtnam and Piazza  \cite{LP2005} and Moriyoshi and Natsume \cite{MN1996}   gave interpretations of the value of the Godbillon-Vey class  in terms of the spectral flow of leafwise Dirac operators for smooth foliations.

The   problem considered in this work was first posed in papers of  Moussu and Pelletier \cite{MoussuPelletier1974} and  Sullivan  \cite{Schweitzer1978}, where they conjectured that a foliation $\F$ with  $GV(\F) \ne 0$ must  have leaves of exponential growth.  The   support  for this conjecture at that time was principally a collection of examples,    and some developing intuition for the dynamical properties of foliations.  
The geometry of  the helical wobble phenomenon is  related  to   geometric properties of contact  flows, such as for the geodesic flow of a compact surface with negative curvature.   The weak stable foliations for such flows have all leaves of exponential growth, and often have non-zero Godbillon-Vey classes \cite{Thurston1972,PlanteThurston1972,ReinhartWood1973, HurderKatok1991,GomesRuggiero2007}.    Moreover, the work of Thurston     in \cite{Thurston1972} implies that for any positive real number $\alpha$   there exist a $C^2$-foliation of codimension-one on a compact oriented $3$-manifold, whose Godbillon-Vey class is $\alpha$ times the top dimension integral cohomology class. These various results suggest that a geometric interpretation of $GV(\F)$ might be related to    dynamical invariants such as ``entropy'', whose values are not limited to a discrete subset of $\mR$.

Given a choice of a complete, relatively compact, 1-dimensional  transversal $\fX \subset M$ to $\F$, the transverse parallel transport along paths in the leaves defines local homeomorphisms of $\fX$, which yields a $1$-dimensional pseudogroup $\GF$ as recalled in Section~\ref{subsec-pseudogroup}. 
The study of the properties of  foliation pseudogroups has been a central theme of foliation theory since the works of Reeb and Haefliger in the 1950's \cite{Reeb1952,Reeb1961,Haefliger1958,Haefliger1962}.

The geometric entropy $\h$ of a $C^1$-foliation $\F$ was introduced by 
Ghys, Langevin and Walczak \cite{GLW1988}, and can be formulated in terms of the pseudogroup $\GF$ associated to the foliation. The geometric entropy is a measure of the  dynamical complexity of the action of $\GF$ on $\fX$, and is one of the most important dynamical invariants of $C^1$-foliations.   The Godbillon-Vey class $GV(\F)$ vanishes for all the known  examples of foliations  for which $\h =0$, and the problem was posed  
  to relate the non-vanishing of the geometric entropy $\h$ of a codimension-one $C^2$-foliation $\F$  with the non-vanishing of its Godbillon-Vey class.

    Duminy   showed in  the unpublished papers \cite{Duminy1982a,Duminy1982b}  that for a $C^2$-foliation of codimension one, $GV(\F) \ne 0$ implies there are leaves of exponential growth. (See   the   account of Duminy's  results in Cantwell and Conlon    \cite{CantwellConlon1984}, and \cite[Theorem~13.3.1]{CandelConlon2003}.) 
Duminy's  proof began by   {assuming} that a $C^2$-foliation $\F$ has no resilient leaves, or equivalently resilient orbits for $\GF$ as in Definition~\ref{defn-resilient}. Then by the 
Poincar\'e-Bendixson theory for codimension--one, $C^2$-foliations \cite{CantwellConlon1984,Hector1983}, Duminy showed  that the Godbillon-Vey class of $\F$ must vanish. Thus, if $GV(\F) \ne 0$  then $\F$ must have at least one resilient leaf.  
If  a codimension--one foliation has a resilient leaf, then by an easy argument it follows that $\F$ has an uncountable set of leaves with exponential growth. Duminy's  proof is ``non-constructive'' and does not directly show how a non-trivial value of the Godbillon-Vey class results in resilient leaves for the foliation.  One of the points of this present work is to give a direct demonstration of this relation, which we show   using techniques of  ergodic theory for $C^1$-foliations.

In the work \cite{GLW1988},  The\'eor\`eme~6.1 states that for a codimension-one, $C^2$-foliation $\F$, if $\h \ne 0$ then $\F$ must have a resilient leaf. Candel and Conlon gave a proof of this result in \cite[Theorem~13.5.3]{CandelConlon2000} for the special case where the foliation is the suspension of a group action on a circle, but were unable to extend the proof to the general case asserted in \cite{GLW1988}. Combining these results, one concludes that for a $C^2$-foliation $\F$, if the geometric entropy $\h =0$ then  $\F$ has no resilient leaves, and thus $GV(\F) = 0$. This result suggests the problem of giving a direct proof of this conclusion.

The development of an ergodic theory approach to the study of the secondary classes  began with the work of Heitsch and Hurder \cite{HeitschHurder1984}, which was inspired by Duminy's  work    \cite{Duminy1982a,Duminy1982b}.
A key idea   introduced in   \cite{Hurder1986,HurderKatok1987}, was to   use techniques from the Oseledets theory of cocycles to study the relation between foliation     dynamics, and the values of  the secondary classes of foliations. 

In this paper, we use     methods from the ergodic  theory of $C^1$-foliations to    show   that for a $C^2$-foliation $\F$, the assumption $GV(\F) \ne 0$ implies that the  foliation $\F$ has resilient leaves, and thus $\h \ne 0$. An important aspect of our proof, is that the subtle techniques of the Poincar\'e-Bendixson theory of $C^2$-foliations are avoided, and the conclusion that there exists resilient leaves follows from straightforward    techniques of dynamical systems.

 The work of Duminy \cite{Duminy1982a} reformulated the study of the Godbillon-Vey class for $C^2$-foliations in terms of the ``Godbillon measure'', which for a $C^1$-foliation $\F$ of a compact manifold $M$, is a linear functional defined on the Borel $\sigma$-algebra $\cB(\F)$ formed from the leaf-saturated Borel subsets of $M$, and by extension this measure is defined on the saturated measurable subsets of $M$. These ideas are introduced and discussed in the papers  \cite{CantwellConlon1984,Duminy1982a,Duminy1982b,HeitschHurder1984,Hurder1986,HurderKatok1987}, and   recalled in  Section~\ref{sec-Gmeasure} below.    Here is our main result, as    formulated in these terms:
\begin{thm}\label{thm1}
If $\F$ is a  codimension--one,  $C^{1}$-foliation  with non-trivial Godbillon measure $G_{\F}$, then $\F$ has a hyperbolic resilient leaf. 
\end{thm}
In the course of our proof of this result,   resilient orbits of the action of the pseudogroup $\GF$ are explicitly constructed using a version of the Ping-Pong Lemma, first introduced by Klein in his study of  subgroups of Kleinian groups \cite{delaHarpe2000}, and which is discussed in Section~\ref{subsec-resilient}.

For $C^{2}$-foliations,  the  Godbillon-Vey class is obtained by evaluating the Godbillon measure on the ``Vey class'' $[v(\F)] | E$ localized to a set $E \in \cB(\F)$. Only the definition of the  class $[v(\F)] | E$  requires that $\F$ be $C^2$.  Thus,  for a $C^2$-foliation $\F$,     $GV(\F) \not= 0$ implies that  $G_{\F} \not= 0$, and we deduce:
\begin{cor}\label{cor1}
If $\F$ is a  codimension--one,  $C^2$-foliation   with non-trivial Godbillon-Vey class $GV(\F) \in H^3(M; \mR)$, then $\F$ has a hyperbolic resilient leaf, and thus the entropy $\h > 0$. 
\end{cor}

We next discuss the strategy of the proof of  Theorem~\ref{thm1}. A key idea in dynamical systems of flows is to consider the points for which the dynamics is ``infinitesimally exponentially expansive'' over long orbit segments, which corresponds to points with positive Lyapunov exponent \cite{BarreiraPesin2013,BDV2005,Pesin1977}. The analog for pseudogroup dynamics  is to introduce the set of points in the transversal $\cX$ for which there are arbitrarily long words in $\GF$  for which the norm of their transverse derivative matrix is exponentially growing   with respect to the word norm on the pseudogroup. 

We introduce in Section~\ref{sec-pesin}, the $\F$-saturated set $\E$ of points in $M$ where the  transverse   derivative cocycle for $\F$ has positive exponent. 
A point $x \in \E \cap \cX$ if and only if there is a sequence of holonomy maps such that the norms of their derivatives at $x$ grow  exponentially fast as a function of ``word length'' in the foliation pseudogroup, and $\E$ is the leaf saturation of this set.

The set $\E$ is a fundamental construction for a $C^1$-foliation.
For example,  a key step in the proof of the generalized Moussu--Pelletier--Sullivan conjecture in \cite{Hurder1986} was to show that for a foliation $\F$ with almost all leaves of subexponential growth,   the Lebesgue measure $|\E| = 0$.   
Here, we   show in Theorem~\ref{thm3} that if  a measurable, $\F$-saturated subset $B \subset M$  is disjoint from $\E$, then the Godbillon measure must vanish on $B$.    

The second step in the proof of Theorem~\ref{thm1}  is to show  that for each point  $x \in \E$,   the holonomy of $\F$ has a uniform exponential estimate along the orbit of $x$ for  its    transverse expansion along arbitrarily long words in the holonomy pseudogroup. This follows from Proposition~\ref{prop-localhol}, which  is pseudogroup version of what is called the ``Pliss Lemma'' in the    literature for non-uniform dynamics \cite{Pliss1972, Mane1987,BDV2005}.   If $\E$ has positive measure,  it is then straightforward to construct  resilient orbits  for the action of $\GF$ on $\cX$, as done    in the proof of    Proposition~\ref{prop-resilient}. The proof of Theorem~\ref{thm1} then follows by combining Theorem~\ref{thm3}, Proposition~\ref{prop-fp} and Proposition~\ref{prop-resilient}.

The proofs of Propositions~\ref{prop-localhol} and \ref{prop-fp} are   the most technical aspects of this paper.
One important issue that arises in the study of pseudogroup dynamical systems,    is that the domain of a holonomy   map in the pseudogroup may depend  upon the ``length'' of the leafwise path   used to define it, so that composing maps in the pseudogroup often results in a contraction of  the domain of definition for the resulting map. 
This is a key difference between the study of dynamics of a group acting on the circle, and that of a pseudogroup associated to a general  codimension--one foliation. One of the key steps in the proof of Proposition~\ref{prop-fp} is  to  show uniform  estimates on the length of the domains of compositions. The proof uses these estimates   to produce an abundance of   holonomy pseudogroup  maps with hyperbolic fixed--points.

We point out one application of Proposition~\ref{prop-fp}, which complements   the main result of \cite{Hurder1991}.

 \begin{thm}\label{thm-noholonomy}
 Let $\F$ be a $C^1$-foliation of codimension-one such that no leaf of $\F$ has a closed loop with hyperbolic transverse  holonomy, then  the hyperbolic set $\E$ is empty.
 \end{thm}

Finally, the extension of the methods for closed manifolds to the case of open manifolds    requires only a minor modification in the definition of the Godbillon measure, as discussed   in Section~\ref{sec-open}.

 For codimension--one foliations, it is elementary that the existence of a resilient leaf implies $\h > 0$.  The converse, that $\h > 0$ implies there is a resilient leaf, was proved in  
 \cite{GLW1988} for $C^2$-foliations, and proved in \cite{Hurder2000c} for $C^1$-foliations. Let ${\rm ``HRL(\F)"} $ denote the property that $\F$ has a hyperbolic resilient leaf. Let $|E|$ denote the Lebesgue measure of a measurable subset $E \subset M$. 
 The results of this paper are     summarized by the following    implications:
\begin{thm}\label{thm-implications}
Let  $\F$ be a  codimension--one,  $C^1$-foliation   of a   manifold $M$. Then
\begin{equation}\label{eq-implications}
 \g \not= 0  \Longrightarrow |\E| > 0 \Longrightarrow {\rm ``HRL(\F)"}  \Longleftrightarrow \h > 0 ~ . 
\end{equation}
\end{thm} 

\medskip

The collaboration of the authors in Spring 1999 leading to this work was made possible by the support of the first author by the Universit\'e of Bourgogne, Dijon. This  support is gratefully acknowledged.

This manuscript is a revised version of a preprint dated March 27, 2004 and submitted for publication \cite{HurderLangevin2004}.  The statements  of the results, and the ideas for their proofs, have not changed in the intervening period, but the   revised manuscript reorganizes the proofs   in Sections~\ref{sec-uniform}  and \ref{sec-hyperbolic},  includes updated references, and   incorporates the suggestions and edits given by the referee of that manuscript.

\section{Foliation Basics}\label{sec-basics}

 In this section, we introduce  some standard notions and results of foliation  geometry and dynamics.
Complete details and further discussions are provided by  the texts  \cite{CN1985,CandelConlon2000,Godbillon1991,HectorHirsch1986,Walczak2004}. 

We   assume that $M$ is a closed oriented smooth Riemannian  $m$-manifold,   $\F$  is a  $C^r$-foliation  of codimension--1 with oriented normal bundle, for $r \geq 1$, and that the leaves of $\F$ are  smoothly immersed submanifolds of dimension $n \geq 2$, where $m = n+1$.  
This is sometimes referred to as a $C^{\infty,r}$-foliation,   where the holonomy transition maps are $C^r$, typically for either $r =1$ or $r=2$.

\subsection{Regular Foliation Atlas}\label{subsec-regular}
 
A   \emph{$C^{\infty,r}$-foliation atlas}  on $M$,  for $r \geq 1$,  is a  finite collection  $\{ (U_{\alpha},\phi_{\alpha})  \mid    \alpha \in \cA \}$ such that:
\begin{enumerate} 
\item $\cU = \{ U_{\alpha}  \mid  \alpha \in \cA  \}$  is an open covering of $M$. 
\item  $\phi_{\alpha} : U_{\alpha} \rightarrow  (-1,1)^m$ is a $C^{\infty,r}$--coordinate chart; that is, for   $(u,w) \in (-1,1)^n \times (-1,1)$, the map $\phi_{\alpha}^{-1}(u,w)$ is $C^{\infty}$ in the ``leaf'' variable $u$ and $C^r$ in the ``transverse'' variable $w$. 
\item Each chart $\phi_{\alpha}$ is transversally oriented.
\item Given  $x \in U_{\alpha} \cap U_{\beta}$ with $\phi_{\alpha}(x) = (u,w)$, for the
  change-of-coordinates map $(u',w') =  \phi_{\beta} \circ \phi_{\alpha}^{-1}(u,w)$, the value of $w'$ is  locally constant with respect to $u$.
\end{enumerate}

\begin{figure}[htbp]
\begin{center}
\includegraphics[width=0.32\textwidth]{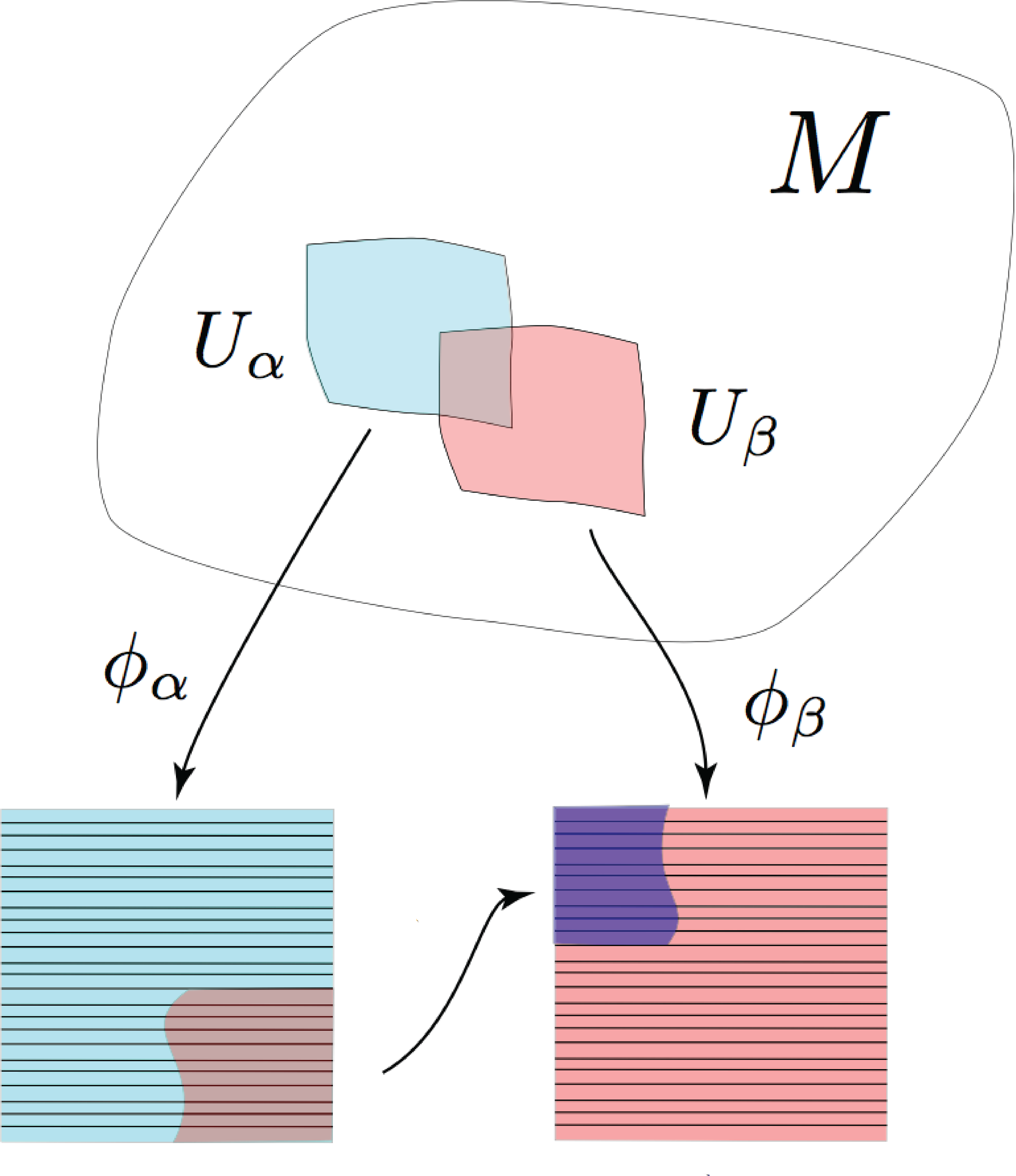}
\caption{Overlapping foliation charts}
\label{Figure1}
\end{center}
\end{figure}

The collection of sets
$$\cV_{\F} \equiv \left \{V_{\alpha, w} =  \phi_{\alpha}^{-1} (V \times \{w\}) \mid V \subset (-1,1)^n ~, ~ w \in (-1,1) ~ , ~ \alpha \in \cA\right\} $$
  form a subbasis for the ``fine topology'' on $M$. 
For $x \in M$, let $L_x \subset M$ denote the connected component of this fine topology containing $x$. Then $L_x$ is path connected, and is called the leaf of $\F$ containing $x$. 
Without loss of generality, we can assume that the coordinates   are positively oriented, mapping the positive orientation for the normal bundle to $T\F$ to the positive orientation on $(-1,1)$.

Note that each leaf $L$ is a smooth, injectively immersed manifold in $M$. The Riemannian metric on $TM$ restricts to a smooth metric on each leaf. The  path-length metric $\dF$ on a leaf $L$ is defined  by   
$$\dF(x,y) = \inf \left\{\| \gamma\| \mid \gamma \colon [0,1] \to L ~ {\rm is } ~ C^1~, ~ \gamma(0) = x ~, ~ \gamma(1) = y  \right\}, $$
  where $\| \gamma \|$ denotes the path length of the $C^1$-curve $\gamma(t)$. If $x,y \in M$   are not on the same leaf, then set $d_{\F}(x,y) = \infty$.
It was    noted by Plante \cite{Plante1975}   that for  each $x \in M$, the leaf $L_x$ containing the point $x$, with the induced Riemannian metric from $TM$   is a complete Riemannian manifold with bounded geometry, that depends continuously on $x$. 
In particular, bounded geometry   implies  that for each $x \in M$, there is a leafwise exponential map
$\exp^{\F}_x \colon T_x\F \to L_x$ which is a surjection, and the composition $\iota \circ \exp^{\F}_x \colon T_x\F \to L_x \subset M$ depends continuously on $x$ in the compact-open topology.

 We  next recall the notion   of a \emph{regular covering}, or what is sometimes called a \emph{nice covering} in the literature  (see \cite[Chapter  1.2]{CandelConlon2000}, or  \cite{HectorHirsch1986}.)
For a regular foliation covering, the  intersections of the coverings of leaves by the plaques of  the charts have nice metric properties. 
 We first recall a standard fact from Riemannian geometry, as it applies to the leaves of $\F$.

For each $x \in M$ and $r > 0$, let $\oB_{\F}(x, r) = \{y \in L_x \mid d_{\F}(x,y) \leq r\}$ denote the closed ball of radius $r$ in the leaf containing $x$.
The Gauss Lemma implies that there exists $\lambda_x > 0$ such that $\oB_{\F}(x, \lambda_x)$ is a \emph{strongly convex} subset for the metric $d_{\F}$. That is, for any pair of points $y,y' \in \oB_{\F}(x, \lambda_x)$ there is a unique shortest geodesic segment in $L_x$ joining $y$ and $y'$ and it is contained in $\oB_{\F}(x, \lambda_x)$ (cf. \cite{BC1964}, \cite[Chapter 3, Proposition 4.2]{doCarmo1992}). Then for all $0 < \lambda < \lambda_x$,  the disk $\oB_{\F}(x, \lambda)$ is also strongly convex.
The compactness of $M$ and the continuous dependence of the Christoffel symbols for a Riemannian metric in the $C^2$-topology on sections of bundles over $M$ yields:
\begin{lemma}\label{lem-stronglyconvex}
There exists $\lF > 0$ such that for all $x \in M$, $\oB_{\F}(x, \lF)$ is strongly convex.
\end{lemma}
If $\F$ is defined by a flow without periodic points, so
that every leaf is diffeomorphic to $\mR$, then the entire leaf is
strongly convex, so $\lF > 0$ can be chosen arbitrarily. For a
foliation with leaves of dimension $n > 1$, the constant $\lF$ must
be less than the injectivity radius for each of the leaves.

Let $d_{M} \colon M \times M \to [0,\infty)$ denote the path-length metric on $M$.  
For $x \in M$ and $\e > 0$, let   $B_{M}(x, \e) = \{ y \in M \mid d_{M}(x, y) < \e\}$ be the open ball of radius $\e$ about $x$, and let 
$\oB_{M}(x, \e) = \{ y \in M \mid d_{M}(x, y) \leq \e\}$ denote its closure.
Then as above,  there exists $\lM > 0$ such that $\oB_M(x, \lambda)$ is a strongly convex ball in $M$ for all $0 < \lambda \leq \lF$.

We use these estimates on the local geometry of $M$ and the leaves of $\F$ to construct a refinement of the given covering of $M$ by foliations charts, which have uniform regularity properties.  

Let $\eU > 0$ be a Lebesgue number for the given covering $\cU$ of $M$.  

Then 
for  each $x \in M$, there exists  $\alpha_x \in \cA$ be such that $x \in B_M(x, \eU) \subset U_{\alpha_x}$.  It follows that 
for each $x \in M$, there exists $0 < \delta_x \leq \lF$ such that $\oB_{\F}(x, \delta_x) \subset B_M(x, \eU)$. 

Let $(u_x, w_x)  = \phi_{\alpha}(x)$, and note  that $\phi_{\alpha}(\oB_{\F}(x, \delta_x)) \subset (-1,1)^n \times \{w_x\}$. 
Then there exists $\e_x > 0$ so that for each $w \in (w_x - \e_x , w_x + \e_x)$ and $x_w = \phi_{\alpha}^{-1}(u_x, w)$ we have 
$\ds \oB_{\F}(x_w, \delta_x) \subset B_M(x, \eU) \subset U_{\alpha_x}$ 
is a leafwise convex subset. Define $U_x$ and $\tU_x $ to be   unions of   leafwise strongly convex disks,
\begin{equation}
U_x   =   \bigcup_{w \in (w_x - \e_x/2 , w_x + \e_x/2)} \, \oB_{\F}(x_w, \delta_x/2) \quad ; \quad 
\tU_x   =   \bigcup_{w \in (w_x - \e_x , w_x + \e_x)} \, \oB_{\F}(x_w, \delta_x)
\end{equation}
  so then $U_x \subset \tU_x \subset B_M(x, \eU) \subset U_{\alpha_x}$. 
The restriction $\phi_{\alpha_x} \colon \tU_x \to (-1,1)^{n+1}$ is then a foliation chart, though the image is not onto. 

Note that for each $x' \in \phi_{\alpha_x}^{-1}(w_x - \e_x , w_x + \e_x)$, the chart $\phi_{\alpha_x}$ defines a framing of the tangent bundle $T_{x'}L_{x'}$  and this framing depends $C^r$ on the parameter $x'$, so we can then use the Gram-Schmidt process to obtain a $C^r$-family of orthonormal frames as well. Then using the inverse of the leafwise exponential map and affine rescaling, we obtain foliation charts 
\begin{eqnarray*}
\wtvp_{\alpha_x} & \colon &  \tU_x \to (-\delta_x,\delta_x)^{n} \times (w_x - \e_x , w_x + \e_x) \cong (-2,2)^n \times (-2,2)\\
\vp_{\alpha_x} & \colon &  U_x \to (-\delta_x/2,\delta_x/2)^{n} \times (w_x - \e_x/2 , w_x + \e_x/2) \cong (-1,1)^n \times (-1,1)
\end{eqnarray*}
where $\vp_{\alpha_x}$ is the restriction of $\wtvp_{\alpha_x}$. Observe that $\wtvp_{\alpha_x}(x) = (\vec{0},0) \in (-1,1)^n \times (-1,1)$ for each $x$.

The collection of open sets  $ \{U_x  \mid x \in M \}$ forms an open cover of the compact space $M$, so there exists a finite subcover ``centered'' at the points $\{x_1, \ldots , x_{\nu}\} \subset M$. Set
\begin{equation}\label{eq-Fdelta}
\dFU = \min \{\delta_{x_1}/2 , \ldots , \delta_{x_{\nu}}/2\} ~ \leq ~ \lF/2 ~.
\end{equation}
This covering by foliation coordinate charts will be fixed and used throughout. 
To simplify notation, for $1 \leq \alpha \leq \nu$, set $U_{\alpha} = U_{x_{\alpha}}$, $\tU_{\alpha} = \tU_{x_{\alpha}}$,    
$\ds  \vp_{\alpha} = \vp_{x_{\alpha}}$,  $\ds  \wtvp_{\alpha} = \wtvp_{x_{\alpha}}$, and   $\cU = \{U_{1}, \ldots , U_{\nu}\}$.

The resulting   collection $\ds \{ \vp_{\alpha}   \colon    U_{\alpha} \to  (-1,1)^n \times (-1,1) \mid 1 \leq \alpha \leq \nu\}$ is a \emph{regular covering} of  $M$ by foliation charts, in the sense used in  \cite[Chapter  1.2]{CandelConlon2000} or  \cite{HectorHirsch1986}.

For each $1 \leq \alpha \leq \nu$, define $\cX_{\alpha} \equiv (-1,1) \cong \{\vec{0}\} \times (-1,1)$ and $\tcX_{\alpha} \equiv (-2,2) \cong \{\vec{0}\} \times (-2,2)$. The extended chart $\wtvp_{\alpha}$ defines     $C^r$--embeddings 
\begin{equation}\label{eq-transversals}
   t_{\alpha} \colon \cX_{\alpha} \rightarrow   U_{\alpha}    \quad , \quad \wtt_{\alpha} \colon \tcX_{\alpha} \rightarrow  \tU_{\alpha} ~ .
\end{equation}
 Let $\fX_{\alpha} = \tau_{\alpha}(\cX_{\alpha})$ and $\tX_{\alpha} = \wtt_{\alpha}(\tcX_{\alpha})$ denote the images of these maps. 
For $n \geq 3$, we can assume without loss of generality that the submanifolds $\tX_{\alpha}$ and $\tX_{\beta}$ are disjoint, for $\alpha \ne \beta$.

Consider   $\cX_{\alpha}$ and $\cX_{\beta}$   as disjoint spaces for $\alpha \ne \beta$, and similarly for $\tcX_{\alpha}$ and $\tcX_{\beta}$. Introduce the disjoint unions of these spaces, as denoted by 
\begin{eqnarray}
 \cX = \bigcup_{1 \leq \alpha \leq \nu} \cX_{\alpha} \quad & \subset  & \quad \tcX = \bigcup_{1 \leq \alpha \leq \nu} \tcX_{\alpha} ~ , \label{eq-transversalT} \\
 \fX = \bigcup_{1 \leq \alpha \leq \nu} \fX_{\alpha} \quad & \subset  & \quad \tX = \bigcup_{1 \leq \alpha \leq \nu} \tX_{\alpha} ~ , \label{eq-transversalX} 
\end{eqnarray}
Note that   $\fX$  is a \emph{complete transversal} for $\F$, as the submanifold $\fX$ is transverse to the leaves of $\F$, and every leaf of $\F$ intersects  $\fX$. The same is true for $\tX$.

 Let $\tau \colon \cX \to \fX \subset M$ denote the map defined by the coordinate chart embeddings $\tau_{\alpha}$, and similarly define $\wtt \colon \tcX \to \tX \subset M$ using the maps $\wtt_{\alpha}$.   

Let each  $\tcX_{\alpha}$ have the metric $\dT$ induced from the Euclidean metric on $\mR$, where   $\dT(x,y) = |x - y|$ for  $x, y \in \tcX_{\alpha}$. Extend this to a metric on $\cX$ by setting $\dT(x,y) = \infty$ for     $x \in \tcX_{\alpha}$,  $y \in \tcX_{\beta}$ with $\alpha \not= \beta$. 
 
 Let each  $\tX_{\alpha}$ have the Riemannian metric  induced from the Riemannian metric on $M$, and let     $\dX$ denote the resulting path-length metric on $\fX_{\alpha}$. As before, extend this to a metric on $\fX$ by setting $\dX(x,y) = \infty$ for     $x \in \tX_{\alpha}$,  $y \in \tX_{\beta}$ with $\alpha \not= \beta$.

Given $r > 0 $ and $x \in \tcX_{\alpha}$ let 
$\ds {\bf B}_{\tX}(x,r)   =  \{ y \in \tX_{\alpha}  \mid   \dX(x,y) < r \}$. Introduce a notation which will be convenient for later work. 
Given a point $x \in \tX_{\alpha}$ and $\delta_1, \delta_2 > 0$,  let 
$$ [x-\delta_1, x +\delta_2] \subset \tX_{\alpha} $$
be the connected closed subset bounded below by the point $x-\delta_1$ satisfying by $\dX(x,x-\delta_1) = \delta_1$ and $[x-\delta_1,x]$ is an oriented interval in $\fX_{\alpha}$.  The set $[x-\delta_1, x +\delta_2]$ is bounded above by the point $x+\delta_2$ satisfying by $\dX(x,x+\delta_2) = \delta_2$ and $[x, x+\delta_1]$ is an oriented interval in $\fX_{\alpha}$.

  For each $1 \leq \alpha \leq \nu$, let $\pi_{\alpha} \equiv \pi_t \circ \vp_{\alpha} \colon U_{\alpha} \to \cX_{\alpha}$ be the composition of the coordinate map $\vp_{\alpha}$ with the projection $\pi_t \colon \mR^{n+1} = \mR^n \times \mR  \to \mR$. 
    For each $w \in \cX_{\alpha}$,   the preimage
$$\cP_{\alpha}(w)   =   \pi_{\alpha}^{-1} \subset U_{\alpha}$$ 
is called a \emph{plaque} of   the chart $\vp_{\alpha}$. For $x \in U_{\alpha}$ we use the notation    $ \cP_{\alpha}(x) = \cP_{\alpha}(\vp_{\alpha}(x))$  to denote the plaque of the chart $\vp_{\alpha}$ containing $x$. Note that $ \cP_{\alpha}(x)$  is  the  connected     component  of the intersection of the leaf $L_x$ of $\F$ through $x$ with the set $U_{\alpha}$. Then the collection of all plaques for the foliation atlas is     indexed by $\cX$.

 The maps $\wtpi_{\alpha} \equiv \pi_t \circ \wtvp_{\alpha} \colon \tU_{\alpha} \to \tcX_{\alpha}$ are defined analogously, with corresponding  plaques $\wtcP_{\alpha}(w)$. For $x \in \tU_{\alpha}$, the plaque of the chart $\wtvp_{\alpha}$ containing $x$ is denoted by $\tcP_{\alpha}(x)   \subset \tU_{\alpha}$.    

Note that each plaque $\cP_{\alpha}(x)$ is strongly convex in the leafwise metric, so if the intersection of two   plaques  $\{ \cP_{\alpha}(x), \cP_{\beta}(y)\}$  is non-empty, then it is   a strongly convex subset.  In particular, the intersection $\cP_{\alpha}(x) \cap \cP_{\beta}(y)$ is   connected. Thus,  each plaque $\cP_{\alpha}(x)$ intersects either zero or one  plaque   in $U_{\beta}$. The same observations  are also true for the extended plaques $\wtcP_{\alpha}(x)$.

\subsection{Holonomy Pseudogroup $\GF$}\label{subsec-pseudogroup}

A pair of indices $(\alpha, \beta)$ is {\em admissible} if $U_{\alpha} \cap U_{\beta} \not= \emptyset$.  For each admissible pair $(\alpha, \beta)$ define 
\begin{eqnarray} 
\cX_{\ab}  & = &  \{x \in \cX_{\alpha}  \mbox{ such that } {\cP}_{\alpha}(x) \cap U_{\beta} \not= \emptyset\} , \label{eq-chiab} \\
\tcX_{\ab}  & = &  \{x \in \tcX_{\alpha}  \mbox{ such that } {\tcP}_{\alpha}(x) \cap \tU_{\beta} \not= \emptyset\} ~ . \label{eq-wtchiab}
\end{eqnarray}
Then there is a well-defined   transition function ${\hh}_{\ba} \colon \cX_{\ab} \to \cX_{\ba}$, which   for $x \in \cX_{\ab} $ is given by 
$$   \hh_{\ba}(x) = y  \mbox{ where  }  {\cP}_{\alpha}(x) \cap {\cP}_{\beta}(y) \not= \emptyset ~ . $$
Note that $\ds  \hh_{\alpha\alpha} \colon \cX_{\alpha} \to \cX_{\alpha}$ is the identity map for each $\alpha \in \cA$.  

The {\it holonomy pseudogroup} $\GF$ associated to the   regular foliation atlas for $\F$ is the pseudogroup with object space $\cX$, and transformations generated by compositions of the local transformations $\ds \{\hh_{\ba} \mid (\alpha,\beta) \mbox{ admissible} \}$.
 The $C^{\infty,r}$--hypothesis on the coordinate charts  implies that each map   $\hh_{\ba}$ is $C^r$. 
Moreover, the hypothesis (2) on regular foliation charts implies that each $\hh_{\ba}$ admits an extension to a  $C^r$-map 
$\ds  \wh_{\ba} \colon \tcX_{\ab} \to \tcX_{\ab}$ 
defined in a similar fashion. The number of admissible pairs is finite, so there exists a uniform estimate on the sizes of the domains of these extensions. We note the following consequence of these observations.

\begin{lemma} \label{lem-leb}
There exists $\e_0 > 0$ so that for every admissible pair $(\alpha, \beta)$ and $x \in \cX_{\ab}$ then   $[x-\e_0, x+ \e_0]  \subset \tcX_{\ab}$. That is, if $x \in \cX_{\alpha}$ is in the domain of $\hh_{\ba}$ then  $[x-\e_0, x+ \e_0]$ is in the domain of $\wh_{\ba}$.
 \end{lemma}
 
For $0 < \delta <  \e_0$  we introduce the closed subsets of $\tcX$ 
\begin{eqnarray}
\cX[\delta]  &   =   &  \{ y \in \tcX \mid \exists ~ x \in \overline{\cX} , ~ \dT(x,y) \leq \delta\} \label{eq-deltanbhd} \\
\cX_{\ab}[\delta]  &   =   &  \{ y \in \tcX_{\ab} \mid \exists ~ x \in \overline{\cX_{\ab}} , ~ \dT(x,y) \leq \delta\} ~ . \label{eq-abdeltanbhd}
\end{eqnarray}
 Thus,  the maps  ${\hh_{\ba}}$ are uniformly $C^r$  on $\cX_{\ab}[\delta]$ for $\delta <  \e_0$.

Composition of elements in $\GF$ will be defined via ``plaque chains''. 
Given   $x,y \in \cX$ corresponding to points on the same leaf, a \emph{plaque chain $\cP$ of length $k$}   between $x$ and $y$ is a collection of plaques
$$ \cP = \{\cP_{\alpha_0}(x_0), \ldots, \cP_{\alpha_k}(x_k) \}, $$ 
where $x_0 = x$, $x_k = y$ and for each $0 \leq i < k$ we have 
$\ds \cP_{\alpha_i}(x_i) \cap \cP_{\alpha_{i+1}}(x_{i+1}) \not= \emptyset$. 
We   write $\| \cP \| = k$. 

A plaque chain $\cP$ also defines an ``extended'' plaque chain for the charts $\{ (\tU_{\alpha}, \widetilde{\phi}_{\alpha})\}$, 
$$\tcP = \{\tcP_{\alpha_1}(x_0), \ldots, \tcP_{\alpha_k}(x_k) \} ~. $$
We say two plaque chains 
$$ \cP  =   \{\cP_{\alpha_0}(x_0), \ldots, \cP_{\alpha_k}(x_k) \} \mbox{  and }  \cQ   =  \{\cP_{\beta_0}(y_0), \ldots, \cP_{\beta_{\ell}}(y_{\ell}) \} $$
are {\it composable} if   $x_k = y_0$, hence $\alpha_k = \beta_0$ and $\cP_{\alpha_k}(x_k) = \cP_{\beta_0}(y_0))$. Their composition is defined by   
$$
\cQ \circ \cP =  \{\cP_{\alpha_0}(x_0), \ldots, \cP_{\alpha_k}(x_k), \cP_{\beta_1}(y_1), \ldots, \cP_{\beta_{\ell}}(y_{\ell})\} ~ .
$$

The holonomy transformation defined by a plaque chain is the local diffeomorphism
\[ \hh_{\cP}  =   \hh_{\alpha_k \alpha_{k-1}} \circ \cdots \circ   \hh_{\alpha_1 \alpha_0} \]
whose domain $\cD_{\cP} \subset \cX_{\alpha_0}$ contains $x_0$. Note that  $\cD_{\cP} $ is the largest connected open subset of $\cX_{\alpha_0}$ containing $x_0$ on which 
$\ds \hh_{\alpha_{\ell} \alpha_{\ell-1}} \circ \cdots \circ \hh_{\alpha_1 \alpha_0}$
is defined for all $0 < \ell \leq k$. The dependence of the domain of $\hh_{\cP}$ on the plaque chain $\cP$ is a subtle issue, yet is at the heart of the technical difficulties arising in the study of foliation pseudogroups. 

Let $\wh_{\tcP}$ be the holonomy associated to the chain 
$ \tcP$, with domain  $\tcD_{\tcP} \subset \tcX_{\alpha_0}$ 
the largest maximal open subset containing $x_0$ on which 
$\ds \wh_{\alpha_{\ell} \alpha_{\ell - 1}} \circ \cdots \circ \wh_{\alpha_1 \alpha_0}$ 
is defined for all $1 < \ell \leq k$.
By the extension property of a regular atlas, the closure $\overline {\cD_{\cP}} \subset \tcD_{\tcP}$ and $\wh_{\tcP}$ is an extension of  $\hh_{\cP}$.

Given a plaque chain
$\ds \cP = \{\cP_{\alpha_0}(x_0), \ldots, \cP_{\alpha_k}(x_k) \}$
and a point $\ds y \in \cD_{\cP}$,  there is 
a ``parallel'' plaque chain denoted 
$\ds \cP(y) = \{\cP_{\alpha_0}(y), \ldots, \cP_{\alpha_k}(y_k) \}$
where $\hh_{\cP}(y) = y_k$.

 For $x \in \cX$, let 
 $$\GF(x) = \{y =  \hh_{\cP}(x) \in \cX \mid ~ \cP ~ {\rm a ~ plaque ~ chain~ for~ which} ~x \in \cD_{\cP}\}$$ 
 denote the orbit of $x$ under the action of the pseudogroup. If $L_{\xi} \subset M$ denotes the leaf containing $\xi \in U_{\alpha}$ with $\pi_{\alpha}(\xi) =x \in \cX_{\alpha}$, then $\tau(\GF(x)) = L_{\xi} \cap \fX$.

\subsection{The derivative cocycle}

Given a plaque chain $\ds \cP = \{\cP_{\alpha_0}(x_0), \ldots, \cP_{\alpha_k}(x_k) \}$ from $x = x_0$ to $y = x_k$,   the  derivative $\hh'_{\cP}(x)$ is defined using the identifications $\cX_{\alpha} = (-1,1)$ for $1 \leq \alpha \leq \nu$. Note that the assumption that the foliation charts are transversally orientation preserving   implies that  $\hh'_{\cP}(x) > 0$ for all plaque chains $\cP$ and $x \in \cD_{\cP}$.

  Given composable plaque chains $\cP$ and $\cQ$, with $x = x_0, y = x_k = y_0, z = y_{\ell}$     the chain rule   implies
\begin{equation}
   \hh'_{\cQ \circ \cP}(x)  =  \hh'_{\cQ}(y) \cdot  \hh'_{\cP}(x) ~ .
\label{eq-cocycle}
\end{equation}
Define the map $\ds  D\hh \colon  \GF \to \mR$   by 
$\ds   D\hh(\cP,y) = \hh'_{\cP(y)}(y)$, which  is called the \emph{derivative cocycle}  for the foliation pseudogroup $\GF$ acting on $\cX$.
The function $\ln \{D\hh(\cP,y) \} \colon  \GF \to \mR$ is   the \emph{additive derivative  cocycle}, or sometimes the \emph{modular cocycle} for $\GF$.

\subsection{Resilient Leaves and Ping-Pong Games}\label{subsec-resilient}

A plaque chain $\ds \cP = \{\cP_{\alpha_0}(x_0), \ldots, \cP_{\alpha_k}(x_k) \}$
is {\it closed} if $x_0 = x_k$. 
A closed plaque chain $\cP$ defines a local diffeomorphism 
$\ds \hh_{\cP} \colon \cD_{\cP} \to \cX_{\alpha_0}$
with $\ds \hh_{\cP}(x) = x$, where $x = x_0 \in \cX_{\alpha_0}$.

A point $y \in  \cD_{\cP}$ is said to be \emph{asymptotic} by iterates of $\hh_{\cP}$ 
to $x$,  if $\ds \hh_{\cP}^{\ell}(y) \in \cD_{\cP}$ for all $\ell > 0$ (where $\ds \hh_{\cP}^{\ell}$ denotes the composition of  $\ds \hh_{\cP}$ with itself $\ell$ times), and  
$\ds   \lim_{\ell \to \infty} \hh_{\cP}^{\ell}(y) = x $. 

The map $\ds \hh_{\cP}$ is said to be a \emph{contraction} at $x$ if there is some $\delta > 0$ so that  every $\ds y \in \bB_{\cX}(x,\delta)$ is asymptotic to $x$.
The map $\ds \hh_{\cP}$ is said to be a \emph{hyperbolic contraction} at $x$ if  $0 < \hh'_{\cP}(x) < 1$. In this case, there exists $\e > 0$ and $0 < \lambda < 1$ so that 
$\ds  \hh'_{\cP}(y) < \lambda$  for all $ y \in  \bB_{\cX}(x,\e)$.
Hence, every point of $\bB_{\cX}(x,\e)$ is asymptotic to $x$, and  
there exists $0 < \delta < \e$ so that  the image of the closed $\delta$--ball about $x$ satisfies
$$ \hh_{\cP}(\overline {\bB_{\cX}(x,\delta)}) \subset  \bB_{\cX}(x,\delta) ~ . $$

\begin{defn}\label{defn-resilient}
We say  $x \in \cX$ is a \emph{hyperbolic resilient point} for $\GF$ if there exists 
\begin{enumerate} 
\item a closed plaque chain $\ds \cP$ such that $\ds \hh_{\cP}$ is  a hyperbolic contraction at $x = x_0$ 
\item a point $y \in  \cD_{\cP}$ which is asymptotic to $x$ (and $y \not= x$)
\item  a plaque chain $\cR$ from $x$ to $y$.
\end{enumerate}
\end{defn}

Figure~\ref{Figure2} below illustrates this concept, where the closed plaque chain $\ds \cP$ is represented by a path which defines it, and likewise for the    plaque chain $\cR$ from $x$ to $y$. Note that the terminal point $y$ is contained in the domain of the contraction $\ds \hh_{\cP}$ defined by $\cP$.

\begin{figure}[htbp]
\begin{center}
\includegraphics[width=0.7\textwidth]{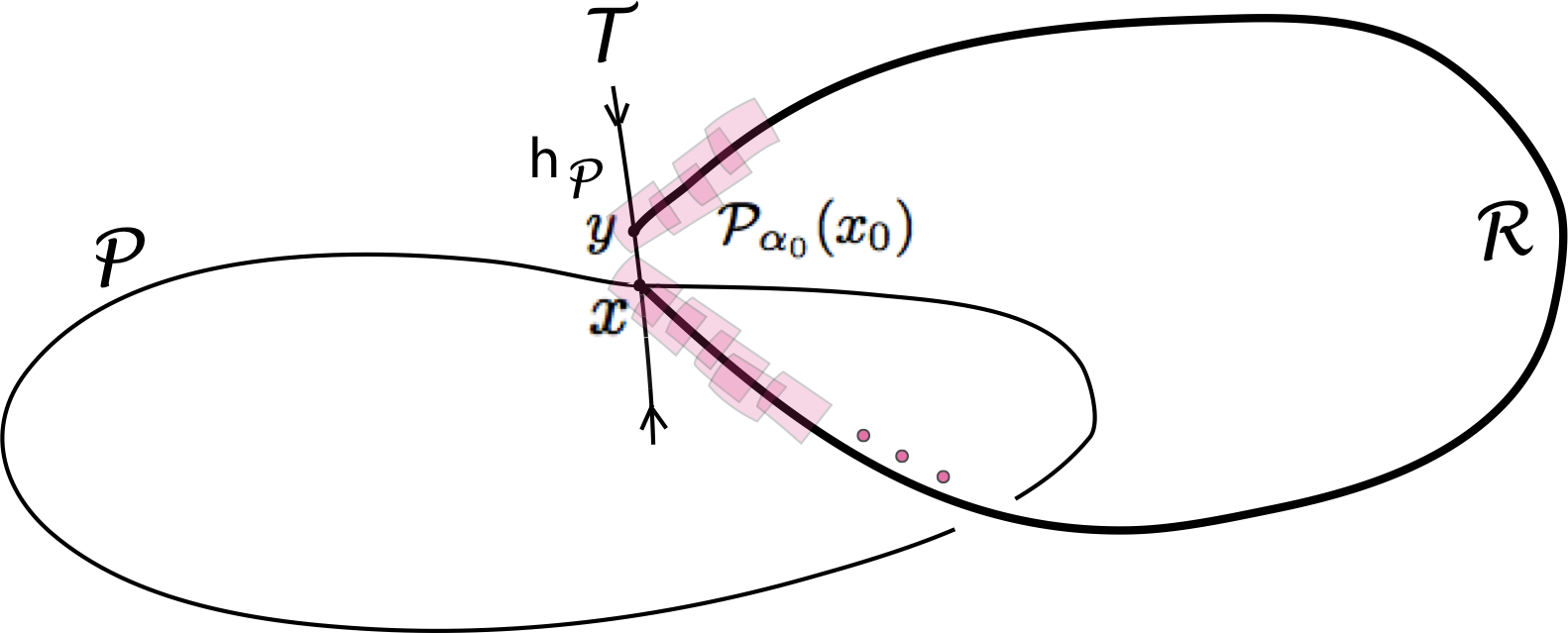}
 \caption{Resilient leaf with     contracting holonomy along loop $\cP$}
\label{Figure2}
\end{center}
\end{figure}

The  ``ping-pong lemma'' is a key technique for the study of 1-dimensional dynamics, which was used by Klein in his study of  subgroups of Kleinian groups \cite{delaHarpe2000}. For a pseudogroup, this has the form:
\begin{defn}
The action of the groupoid  $\GF$  on $\cX$ has a   ``\emph{ping-pong game}'' if  there exists $x,y \in \cX_{\alpha}$ with $x \not= y$ and
\begin{enumerate} 
\item a closed plaque chain $\ds \cP$ such that $\ds \hh_{\cP}$ is  a   contraction at $x = x_0$ 
\item  a closed plaque chain $\ds \cQ$ such that $\ds \hh_{\cQ}$ is  a   contraction at $y = y_0$ 
\item $y \in \cD_{\cP}$ is asymptotic to $x$ by $\ds \hh_{\cP}$ and $x\in \cD_{\cQ}$ is asymptotic to $y$ by $\ds \hh_{\cQ}$
\end{enumerate}
\label{defn-ppg}
We say that the ping-pong game is {\em hyperbolic} if the   maps  
$\ds \hh_{\cP}$ and $\ds \hh_{\cQ}$ are hyperbolic contractions.
\end{defn}

Figure~\ref{Figure3} below illustrates  the ping-pong dynamics, where the closed plaque chain $\ds \cP$ is represented by a path which defines it, and likewise for the    plaque chain $\cQ$.  

\begin{figure}[htbp]
\begin{center}
\includegraphics[width=0.7\textwidth]{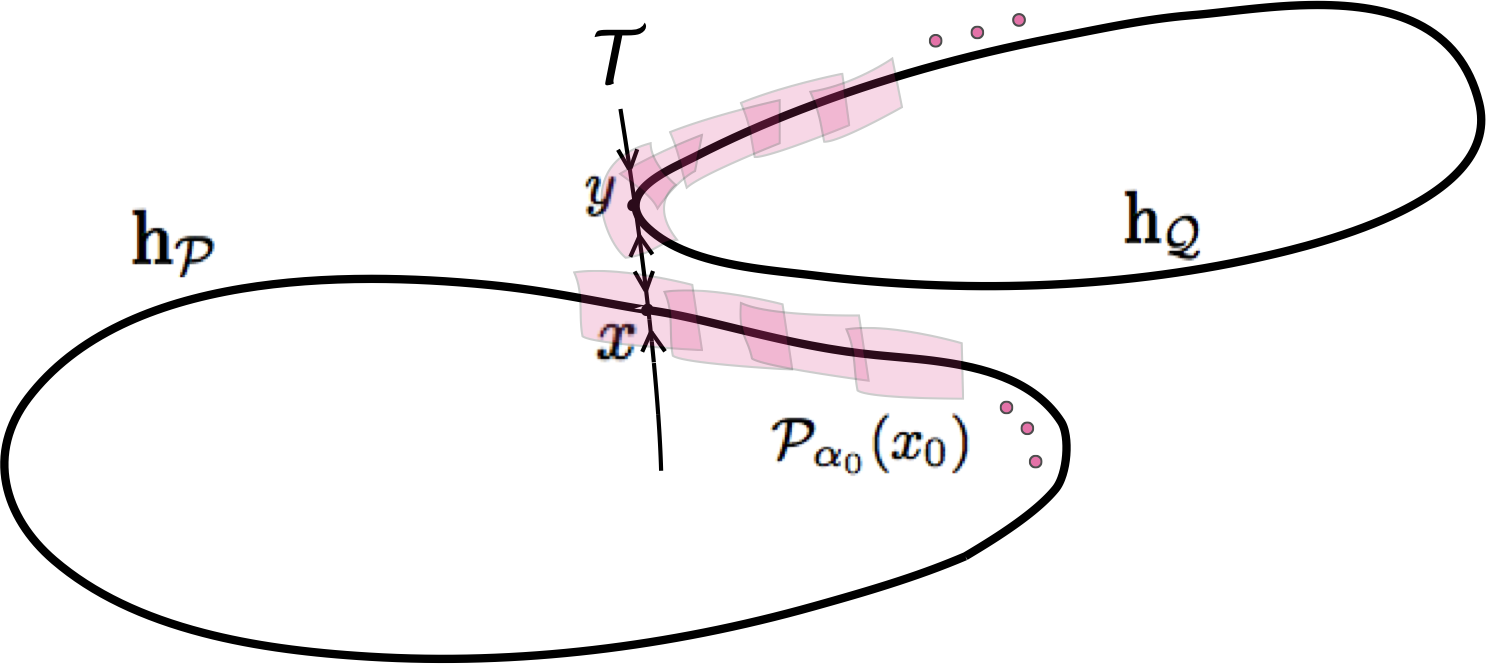}
\caption{Closed paths $\cP$ and $\cQ$ with contracting holonomy   generate  a ping-pong game}
\label{Figure3}
\end{center}
\end{figure}

These two notions are closely related as follows; for example, see   \cite{GLW1988} for a more detailed discussion.
\begin{prop} \label{prop-ppg1}
 $\GF$  has a ``ping-pong game'' if and only if  it has a   resilient point, and has 
a ``hyperbolic ping-pong game'' if and only if  it has a  hyperbolic resilient point.
\end{prop}

\section{The Godbillon-Vey invariant}\label{sec-Gmeasure}

In this section, we recall the definition of the Godbillon-Vey class, the basic ideas of the Godbillon measure, and how it is used to estimate the values of the Godbillon-Vey invariant. The textbook by Candel and Conlon \cite[Chapter 7]{CandelConlon2003} gives a detailed discussion   of the Godbillon-Vey class and its properties. 
Proposition~\ref{prop-vanishing} is the key result of this section that we will use to relate the Godbillon-Vey invariant to the dynamics of the foliation. 

\subsection{The Godbillon-Vey class}
The  Godbillon-Vey class is well-defined for $C^2$-foliations, and the Godbillon measure for $C^1$-foliations. However, giving these definitions for $C^{r}$-foliations,  for $r = 1$ or $2$, adds a layer of notational complexity which obscures the basic ideas of  the constructions. Thus, for clarity of the exposition, 
we assume throughout \emph{this section} that $\F$ is a $C^{\infty}$-foliation, and leave to the reader the technical modifications to the arguments which are required to show the analogous results for $C^r$-foliations, or the reader may consults the works \cite{Duminy1982a,HeitschHurder1984,Hurder1986}. 

Assume that $M$ has a Riemannian metric, and that $\F$ is a $C^{\infty}$-foliation of codimension one,  so   that $M$ admits a covering by smooth regular foliation charts. Thus, the normal bundle $Q \to M$ to $T\F$ can be identified with  the orthogonal space to the tangential distribution $T\F$.  We may assume without loss of generality that $M$ is connected, and that both the tangent bundle $TM$ and the normal bundle $Q$ are oriented, as the dynamical properties to be studied are preserved after passing to a finite covering of $M$.  Thus,  $T\F$ is defined as the kernel of a non-vanishing $1$-form $\omega$ on $M$. 

As  the distribution $T\F$ is integrable, the Froebenius Theorem implies that $d\omega \wedge \omega = 0$. This property is used to define the Godbillon-Vey class, as described below. It is also well-known from the foliation literature that integrability   allows for the construction of various differential graded algebras which are  derived from the de~Rham complex $\Omega^*(M)$ of $M$, each of which  which reflect aspects of the geometry of $\F$. These complexes and their properties are fundamental for the definition of the Godbillon measure, following the works \cite{Duminy1982a,HeitschHurder1984}.

 We first recall a basic construction that is used throughout the following discussions. 
 Let $\omega$ be a non-vanishing $1$-form on $M$ whose kernel equals $T\F$, and $\vec{v}$ a vector field on $M$ such that $\omega(\vec{v}) =1$.
  The integrability of the tangential distribution $T\F$ implies that  $d\omega \wedge \omega = 0$. 
Hence,  there exists a $1$-form $\alpha$ with $d\omega =  \omega \wedge  \alpha$. 
The choice of the $1$-form $\alpha$ is not canonical,  and introduce a procedure for choosing a representative   for $\alpha$.  Set   
$\eta = \iota(\vec{v}) d\omega$, and note  that $\eta(\vec{v}) = 0$.
Then for any  choice of $\alpha$ such that $d\omega =  \omega \wedge  \alpha$, let        $\vec{u}$ be tangent to $\F$,   then we have 
\begin{equation}\label{eq-canonical}
\eta(\vec{u}) =  (\iota(\vec{v}) d\omega) (\vec{u}) = d\omega(\vec{v}, \vec{u})  = (\omega \wedge \alpha) ( \vec{v}, \vec{u})     = \alpha(\vec{u})
\end{equation}
 as $\omega(\vec{v}) = 1$ and $\omega(\vec{u}) = 0$ by definition. Thus, for any $1$-form $\alpha$ such that $d\omega =  \omega \wedge  \alpha$ and any leaf $L$ of $\F$, we have that their restrictions satisfy $\alpha|L = \eta |L$. 
 We introduce    the following   notation.
  
 \begin{defn} \label{def-dF} Let $\omega$ be a non-vanishing $1$-form on $M$ whose kernel equals $T\F$, and $\vec{v}$ a vector field on $M$ such that $\omega(\vec{v}) =1$. Define $\Dv \omega =   \iota(\vec{v}) \, d\omega$.   
 \end{defn}

Then with this definition, we have  $d\omega = \omega \wedge \Dv \omega$.  
 The restricted leafwise $1$-form $\Dv \omega |\F \colon T\F \to \mR$ is known in the literature as the \emph{Reeb form}, and first appeared in the works of Reeb  \cite{Reeb1952,Reeb1961}, and was even implicitly introduced by Poincar\'e \cite{Poincare1881}. The leafwise form $\Dv \omega|\F$ has an interpretation as  the gradient of the \emph{Radon-Nikod\'yn derivative}  along leaves for the ``transverse measure'' for $\F$ defined by the $1$-form $\omega$, as discussed in  
  \cite{Ghys1989,HeitschHurder1984,Hurder1986}. Furthermore,  this is the idea behind  the relation between the dynamics of $\F$ with the flow-of-weights for the \emph{von Neumann algebra} associated to $\F$, as discussed by Connes in \cite{Connes1994}.

Now set    $\eta = \Dv \omega$, and calculate
\begin{equation}\label{eq-etaideal}
0 = d(d\omega) = d( \omega \wedge \eta) = d\omega \wedge \eta  - \omega \wedge d\eta  = \omega \wedge \eta  \wedge \eta  - \omega \wedge d\eta  = - \omega \wedge d\eta ~. 
\end{equation}
shows that  $d\eta \wedge \omega =0$, and hence the $2$-form $d\eta$ is a multiple of $\omega$.  Then calculate  $d(\eta \wedge d\eta) = d\eta \wedge d\eta = 0$   as $\omega \wedge \omega = 0$, so  that $\eta \wedge d\eta$ is a closed 3-form.   

Throughout this work, $H^*(M)$ will denote the de~Rham cohomology groups of $M$.

\begin{thm}[Godbillon and Vey, \cite{GodbillonVey1971}]\label{thm-GV}
The cohomology class $GV(\F) = [\eta \wedge d\eta] \in H^3(M)$ is independent of the choice of the $1$-forms $\omega$ and $\eta$.
\end{thm}

Moreover, the Godbillon-Vey class $GV(\F)$ is an  invariant of the   foliated concordance class of $\F$, as noted for example in  
Thurston \cite{Thurston1972} and Lawson \cite[Chapter~3]{Lawson1977}.  
 
 The definition of the Godbillon-Vey class in Theorem~\ref{thm-GV} reveals very little about the relation of this cohomology class with the dynamics of the foliation $\F$. 
 In the case where the leaves of $\F$ are defined by a smooth fibration $M \to \mS^1$,  the defining $1$-form $\omega$ for $\F$ can be chosen to be a closed form, and  it is then immediate from the definition that   $GV(\F) =0$.  For codimension-one foliations with   slightly more dynamical complexity, the proof that $GV(\F) =0$ 
using the definition above becomes far more involved. For example, 
 Herman showed in \cite{Herman1977} that a foliation defined by the suspension of an action of the abelian group $\mZ^2$ on the circle must have   $GV(\F) =0$. The  proof used an averaging process to obtain a sequence of defining $1$-forms $\omega_n$ for which the corresponding $1$-forms $\Dv \omega_n \to 0$.
 Subsequently,   $GV(\F) = 0$ was shown for the progressively more general  classes of  \emph{foliations without holonomy} by  Morita and Tsuboi \cite{MoritaTsuboi1980},  for foliations \emph{almost without holonomy} by T. Mizutani, S. Morita and T. Tsuboi  \cite{MMT1981}, 
and the   case of foliations which admit \emph{SRH decompositions} by Nishimori \cite{Nishimori1980} and Tsuchiya \cite{Tsuchiya1982}.

The breakthrough idea of Duminy, which can be seen first in his paper with Sergiescu \cite{DuminySergiescu1981}, and further developed in the unpublished work \cite{Duminy1982a}, is to introduce the notion of the  \emph{Godbillon functional}, and its strategy    to separate the roles of the forms $\eta$ and $d\eta$
in the definition of $GV(\F)$, and then the study of how the contribution from the form $\eta$ is related to the dynamical properties of $\F$. 
The definition of the Godbillon functional  requires considering    cohomology classes defined by the forms $\eta$ and  $d\eta$ in their ``largest possible''  natural contexts. 
Introduce the spaces, for $p \geq 1$, 
\begin{equation}
A^p(M,\F) \equiv \{ \xi = \omega \wedge \beta \mid \beta \in \Omega^{p-1}(M) \}  \subset \Omega^*(M) .
\end{equation}
The space    $A^p(M,\F) $ can alternately be defined as the space of $p$-forms on $M$ which vanish when restricted to each leaf of $\F$. Note that the identity $\omega \wedge \omega = 0$, implies that the product of forms in $A^p(M,\F)$ and $A^p(M,\F)$ always vanishes,  so the sum of these spaces is a subalgebra $A^*(M,\F)  \subset \Omega^*(M)$. We next show that this is a differential subalgebra.
  
  The identity
$d\omega = \omega \wedge \eta$ implies that $A^*(M,\F)$ is closed under exterior differentiation. More precisely, 
let  $\xi =  \omega \wedge \beta \in A^k(M,\F) $ for $k \geq 1$, then  
 \begin{equation}\label{eq-leafwisecomplex}
d\xi = d(\omega \wedge \beta) = d\omega \wedge \beta - \omega \wedge d\beta =   (\omega \wedge \eta) \wedge \beta  - \omega \wedge d\beta = \omega \wedge (\eta - d\beta) \in A^{k+1}(M,\F) . 
\end{equation}
Thus,  $A^*(M,\F)$ is a differential graded algebra.  
Let  $H^{*}(M,\F)$ denote the    cohomology   of the differential graded complex $\{A^*(M,\F), d \})$. For a closed form $\xi \in A^k(M,\F)$,  let $[\xi]_{\F} \in H^{k}(M,\F)$ denote its cohomology class.
  
  The  calculation \eqref{eq-leafwisecomplex} shows that differential on the complex $A^k(M,\F) $ is ``twisted'' by the $1$-form $\eta$. Twisted cohomology also arises in the study  of the dynamics of Anosov flows by   Fried in \cite{Fried1986}, and it would be interesting to understand if  (partial)  results analogous to those in \cite{Fried1986}  can be obtained from the study of the  cohomology spaces $H^*(M,\F)$.

The inclusion of the ideal $ A^*(M,\F) \subset  \Omega^*(M)$ induces a map on cohomology 
$H^{*}(M,\F)    \to   H^{*}(M)$. In general, the induced map need not be injective,  and 
the calculation of the cohomology groups $\ds  H^{*}(M,\F)$ is often an intractable problem \cite{ElKacimi1983}. However, it is the fact that   $\ds  H^{*}(M,\F)$ is the domain of the Godbillon operators which makes them important. We next discuss the linear functionals defined on these spaces.

First, we make explicit a property implied by the discussions above. 
Let    $\eta$ be any choice of a $1$-form satisfying $d\omega =  \omega \wedge  \eta$. Let    $\xi \in A^k(M,\F)$ with $k \geq q$ so that  $\xi = \omega \wedge \beta$ for some $(k-q)$-form $\beta \in \Omega^{k-q}(M)$. 
Then by the calculation \eqref{eq-canonical}, the product $\eta \wedge \xi = \eta \wedge \omega \wedge \beta$
    depends only on the leafwise restriction of the form $\eta$. 
Thus,  $\eta \wedge \xi$  is   independent of the choice of $\eta$ which satisfies $d\omega =  \omega \wedge  \eta$, and in particular, it equals  $\Dv \omega \wedge \xi$ where $\vec{v}$ is a vector field on $M$ such that $\omega(\vec{v}) =1$. 

  Again, let    $\eta$ be any choice of a $1$-form satisfying $d\omega =  \omega \wedge  \eta$.
 Recall from \eqref{eq-etaideal} that the closed $2$-form $d\eta$ is in the ideal generated by $\omega$, so $d\eta \in A^2(M,\F)$. 
Duminy observed in  \cite{Duminy1982a} (see also \cite[Chapter 7]{CandelConlon2003},\cite{HeitschHurder1984})  that    the class $[d\eta]_{\F} \in H^{2}(M,\F)$   is independent of the choice of the   $1$-form $\eta$, and so is an invariant of $\F$,  which he called  the {\it Vey class} of $\F$.
 The $2$-form $d\eta$ has some   properties analogous to those of a symplectic form on $M$, especially in the geometric interpretation of the Godbillon-Vey invariant as ``helical wobble'' \cite{LW2008,ReinhartWood1973,Thurston1972}, but the analogy is very loose. The geometric meaning of the   class 
$[d\eta]_{\F}$ remains obscure, although as noted below, $[d\eta]_{\F} = 0$ implies that $GV(\F) = 0$. 

\subsection{The Godbillon operator}
Let a defining $1$-form $\omega$ be given, and a vector field on $M$ such that $\omega(\vec{v}) =1$, and set $\eta = \Dv \omega =   \iota(\vec{v}) d\omega$. 
 Given a  closed form  $\xi \in    A^{p}(M,\F)$,  the product 
$\eta \wedge \xi \in  A^{p+1}(M,\F)$ is   closed, as $d(\eta \wedge \xi) = d\eta \wedge \xi = \omega \wedge  \eta \wedge \omega = 0$.
Moreover, if $\xi = d\beta$ for some form $\beta \in A^{p-1}(M,\F)$,  then 
$\eta \wedge \beta \in A^{p}(M,\F)$ and 
\begin{equation}
d(-\eta \wedge \beta) =  -(d\eta) \wedge \beta + \eta \wedge d\beta = \eta \wedge \xi ~ .
\end{equation}
Thus, given $[\xi]_{\F} \in H^{p+1}(M, \F)$ we obtain a well-defined class  
 $g([\xi]_{\F}) = [\eta \wedge \xi]_{\F} \in H^{p+1}(M, \F)$. 
It follows   that there is a well-defined composition 
\begin{equation}\label{eq-pairing}
g   \colon   H^{p}(M,\F)    \to H^{p+1}(M,\F)    \to   H^{p+1}(M) ~ , ~~ g([\zeta]_{\F}) = [\eta \wedge \zeta] 
\end{equation}
which is called the   \emph{Godbillon operator}.  It was shown above  that the $2$-form $d\eta$ is a multiple of $\omega$, and is clearly a closed form, so it defines a cohomology  class $[d\eta]_{\F} \in H^{2}(M,\F)$. Then we have 
  $g([d\eta]_{\F}) = [\eta \wedge d\eta] = GV(\F) \in H^3(M)$. That is,  ``Godbillon(Vey) = Godbillon-Vey''.

  If $M$ is a closed $3$-manifold with fundamental class $[M]$, then evaluating  $GV(\F)$  on $[M]$  yields a real number,     the \emph{real Godbillon-Vey invariant} of $\F$:
  $$ \langle GV(\F), [M] \rangle = \int_M ~ \eta \wedge d\eta ~ .$$

If $M$ is an open $3$-manifold, then $H^3(M) = 0$ so that $GV(\F) = 0$ in this case. However, the class $GV(\F)$ need not vanish in the case  when $M$ is open and  $M$ has dimension $m > 3$. In this case,   it is necessary to introduce cohomology with compact supports, in order   to obtain a real-valued invariants from   the class $GV(\F)$.

 Now let  $\Omega_c^*(M) \subset \Omega^*(M)$ denote the differential subalgebra of forms with compact support. The cohomology of this ideal is denoted by 
  $H_c^*(M)$ which is called with the \emph{de~Rham cohomology with compact supports} of $M$. Let $A_c^*(M, \F) \subset \Omega_c^*(M)$ denote the differential ideal consisting of forms in $A^*(M, \F)$ with compact support.  Its cohomology groups are denoted by  $H_c^*(M, \F)$, and these groups are called the \emph{foliated cohomology with compact supports}.

Given  a closed form $\zeta  \in A^{p}(M, \F)$,    let $\xi \in   \Omega_c^{k}(M)$ be  a closed form with compact support, 
then  the product $\zeta \wedge \xi \in A^{k+p}(M,\F)$ is again closed with compact support.
If  either form is the boundary of a   form with compact supports, then  $\psi \wedge \xi$ is also the boundary of a compact form. Thus, there is a well-defined pairing 
\begin{equation} \label{pairing}
H^{p}(M,\F) \times H_c^k(M)     \to H_c^{k + p}(M, \F) ~ .
\end{equation}
In particular, given a class $[\xi] \in H_c^{m-3}(M)$ represented by a smooth closed form $\xi \in \Omega_c^{m-3}(M)$,   then   the pairing $[d\eta]_{\F} \cup [\xi] = [d\eta \wedge \xi]_{\F} \in H_c^{m-1}(M,\F)$ is well-defined.

 Recall that the manifold $M$ is assumed to be oriented and connected, so by Poincar\'e duality      the pairing  $H^p(M) \otimes H_c^{m-p}(M) \to H_c^m(M) \cong \mR$ is non-degenerate for $0 \leq p < m$.
In particular, the value of the  class $ [\eta \wedge d\eta] \in H^3(M)$ is determined by its pairings with classes in $H_c^{m-3}(M)$. This is the idea behind the next concept, which is the basis for the Godbillon measure.

The  Godbillon operator  in \eqref{eq-pairing} applied to a class in $H_c^{m-1}(M,\F)$ yields a closed $m$-form with compact support on $M$,  which can be integrated over the fundamental class  to obtain a real number.  This composition yields a linear functional    denoted by 
\begin{equation}\label{pd}
G \colon H_c^{m-1}(M,\F)\to   \mR, \;\; ~~ \;\;  G([\zeta]_{\F}) = \langle [\eta \wedge \zeta], [M] \rangle = \int_M \eta \wedge \zeta ~ .
\end{equation}
Note that we use the notation ``$g$'' for the Godbillon operator between cohomology groups, and the notation ``$G$'' for the linear functional on the cohomology group $H_c^{m-1}(M,\F)$. 

With these preliminary preparations, we have the basic result:

\begin{prop}[Duminy, \cite{Duminy1982a}]\label{prop-duality}
The value of the Godbillon-Vey class $GV(\F) \in H^3(M)$ is determined by the Godbillon operator $G$ in \eqref{pd}. In particular, if $G \equiv 0$ then $GV(\F) = 0$.
\end{prop}
\proof
For the case when the dimension $m =3$ and $M$ is compact, this follows by applying the linear functional $G$ to the class $[d\eta]_{\F} \in H^{2}(M, \F) = H_c^{2}(M, \F)$. 
For  $m > 3$, then by  Poincar\'e duality, the value of $GV(\F) \in H^3(M)$ is determined by  pairing the $3$-form $\eta \wedge d\eta$ with closed forms   $\xi \in \Omega_c^{m-3}(M)$, followed by integration, to obtain 
$$\langle GV(\F) \cup [\xi], [M] \rangle  = \int_M ~ (\eta \wedge d\eta) \wedge \xi  .$$
Note that  $[d\eta \wedge \xi]_{\F} \in H_c^{m-1}(M,\F)$, so  that $\langle GV(\F) \cup [\xi], [M] \rangle = G([d\eta \wedge \xi]_{\F})$.  The claim follows.
\endproof

   This elementary    observation by Duminy    implies    that  $GV(\F) = 0$ if     $G$  is the trivial functional.  The strategy to proving that $GV(\F) = 0$ is thus, to  obtain dynamical properties of a foliation which suffice to show that the   linear functional  $G$ vanishes.

\subsection{The Godbillon measure}
 The linear  functional $G$ possesses special  properties that were hinted at in the literature preceding Duminy's work (see the survey \cite{Hurder2002} for a fuller discussion of the ideas leading up to Duminy's work.)   In particular, Duminy   showed that the integrand in \eqref{pd} which defines the operator $G$ can be restricted to   saturated Borel subsets, to obtain well-defined real invariants. This   observation was systematically generalized in the work        
  \cite{HeitschHurder1984},  to show that    $G$ extends to a generalized measure on the \emph{Lebesgue measurable} saturated subsets of $M$. Moreover, the  values of the measure  can be calculated using \emph{measurable cocycle data}, as discussed in \cite{Hurder1986}.   The extension to measurable data allows the introduction of techniques of ergodic theory. We discuss the definition of the Godbillon measure and its properties in more detail below.

A set $B \subset M$ is $\F$--saturated if for all $x \in B$,  the leaf $L_x$ through $x$ is contained in $B$. Let $\cB(\F)$  denote the $\Sigma$-algebra of Lebesgue measurable $\F$--saturated subsets of $M$.  

\begin{thm} \cite{Duminy1982a, HeitschHurder1984}\label{thm-duminy}
For each $B \in \cB(\F)$, there is a well-defined linear functional 
\begin{equation}\label{eq-measure}
 G_{\F}(B) \colon H_c^{m-1}(M,\F) \to  \mR \quad , \quad G_{\F}(B)([\zeta]_{\F}) = \int_B \eta \wedge \zeta
\end{equation}
where $\zeta \in A_c^{m-1}(M,\F)$ is closed. Note that if $B$ has Lebesgue measure zero, then $G_{\F}(B) =0$.
Moreover, the correspondence 
$$ B \mapsto G_{\F}(B) \in {\rm Hom_{cont}}(H_c^{m-1}(M,\F), \mR)$$
is a countably additive measure on $\cB(\F)$, called the \emph{Godbillon measure}.
\end{thm}

Part of the claim of Theorem~\ref{thm-duminy} is that the linear functional (\ref{eq-measure})
is independent of the choice of the smooth $1$-form $\omega$ defining $\F$. Much more is true, as described below.
The key idea, introduced in \cite{HeitschHurder1984}, is to consider the space of leafwise forms on $\F$ which are leafwise smooth, but need only be measurable as functions on $M$.

 Introduce the graded differential algebra   $\Omega^*(\F)$ consisting of leafwise forms. That is, for $k \geq 0$, the space $\Omega^k(\F)$   consists of sections of the dual to the $k$-$th$ exterior power of the leafwise tangent bundle $T\F$. 
 That is, a form $\xi \in \Omega^k(\F)$ is defined, for each $x \in M$, on a $k$-tuple $(\vec{v}_1, \ldots, \vec{v}_k)$ of vectors in the tangent space $T_x\F$ to the leaf $L_x$ containing $x$. Moreover, we   require that for any leaf $L$ of $\F$, the restriction $\xi |L$  to $L$ is a smooth form.  
There is   a leafwise exterior differential 
\begin{equation}\label{eq-defDF}
\DF \colon \Omega^k(\F) \to \Omega^{k+1}(\F) \quad , \quad \DF(\xi) = d(\xi | L).
\end{equation}
For  $\xi \in \Omega^k(\F)$, the definition of $\DF(\xi) \in \Omega^{k+1}(\F)$ is as follows. For each leaf $L$ of $\F$, 
 the restriction $\xi | L$  is a smooth $k$-form on $L$, so there is a well-defined  exterior differential $d(\xi | L)$. The collection of leafwise forms $\{ d(\xi | L) \in \Omega^{k+1}(L) \mid L \subset M\}$
  defines the class $\DF(\xi) \in \Omega^{k+1}(\F)$.  The cohomology groups of the graded differential algebra $\{\Omega^*(\F), \DF\}$ are called the \emph{foliated cohomology} of $\F$.
  
 A key observation in the definition of the  exterior differential in \eqref{eq-defDF} is that it does not require any regularity for the transverse behavior of the leafwise forms. Thus, one can consider the subcomplex $\Omega_{\infty}^*(\F) \subset \Omega^*(\F)$ of smooth leafwise forms and the corresponding space $H_{\infty}^*(\F)$ of smooth foliated cohomology, which was used by Heitsch in \cite{Heitsch1975} to study the deformation theory of foliations. If the forms are assumed to be continuous, we obtain 
  the subcomplex $\Omega_{c}^*(\F) \subset \Omega^*(\F)$ whose cohomology spaces  $H_{c}^*(\F)$ were  studied by  El Kacimi-Alaoui in \cite{ElKacimi1983}.
One can also consider the subcomplex $\Omega_{m}^*(\F) \subset \Omega^*(\F)$ of  measurable (or bounded measurable) sections of the dual to the $k$-$th$ exterior power of the leafwise tangent bundle $T\F$, then we obtain the  measurable cohomology leafwise cohomology $H_{m}^*(\F)$ groups used by Zimmer in  \cite{Zimmer1981,Zimmer1990}  to study the rigidity theory for measurable group actions.  The next result we require is formulated using the complex subcomplex $\Omega_{m}^*(\F)$.
 
A function $f \colon M \to \mR$ is said to be \emph{transversally measurable} if it is a measurable function, and for each leaf $L$ of $\F$, the restriction $f | L$  is smooth   and the leafwise derivatives of $f$ are measurable functions as well.  Such a function $f$ is the typical element in $\Omega_m^0(\F)$. Given     $f \in \Omega_m^0(\F)$ and a form $\xi \in \Omega_{c}^k(\F)$, then the product $f \cdot \xi \in \Omega_m^k(\F)$. We next   introduce norms on the spaces $\Omega_m^k(\F)$.

For each $x \in M$, the Riemannian metric on $T_x M$ defines a norm on $T_x M$, which restricts to a norm on the leafwise tangent space $T_x\F$. The norm on the space $T_x\F$  induces a dual norm on the cotangent bundle $T_x^* \F$, and also induces norms on each exterior vector space $\Lambda^k T_x \F$ and on its dual $\Omega^k(T_x \F)$, for all $k > 1$. We denote this norm by  $\| \cdot \|_x$ in each of these cases. For a function $f \in \Omega_x(\F)$, let $\| f\|_x = |f(x)|$.

Given a subset $B \subset M$, and a leafwise  form $\xi \in \Omega^k(\F)$ for $k \geq 0$, define the sup-norm over $B$ by   
$$\|\xi\|_{B} = \sup_{x \in B} ~ \|\xi_x\| ~.$$

A remarkable property of the Godbillon measure $G_{\F}$, as shown in Theorem~2.7 of \cite{HeitschHurder1984}, is that for $B \in \cB$ the value of $G_{\F}(B)$   can   be calculated using a $1$-form $\omega_f = \exp(f) \cdot \omega$, where we require that $f \in \Omega_m(\F)$, and for $\vec{v}$ with $\omega_f (\vec{v}) = 1$ and $\eta_f = \Dv (\omega_f)$, we have   $\| \eta_f \|_B < \infty$. 
Then \cite[Theorem~2.7]{HeitschHurder1984} shows that  given a closed form with compact support $\zeta  \in A_c^{m-1}(M, \F)$,   
\begin{equation}\label{eq-integralcalc}
G_{\F}(B)([\zeta]_{\F}) ~ = ~    \int_B ~ \eta_f \wedge \zeta ~ .
\end{equation}
We now recall a fundamental result, Theorem~4.3 of \cite{Hurder1986}, which is a broad generalization of the ideas in the seminal work by Herman \cite{Herman1977}:

\begin{prop} \label{prop-vanishing}
Let  $B \in \cB(\F)$. Suppose there exists a sequence of transversally measurable functions $\{f_n \mid n =1,2, \ldots\}$ on $M$ so that   the    1-forms  $\{\omega_n = \exp(f_n) \cdot \omega \mid n = 1, 2, \ldots\}$ on $M$  satisfy     $\| D^{\vec{v}_n}(\omega_n)\|_B < 1/n$ where $\omega_n(\vec{v}_n) =1$.  Then $G_{\F}(B) = 0$.  
\end{prop}
\proof
For each $n \geq 1$, let $\vec{v}_n$ be a vector field on $M$ such that $\omega_n(\vec{v}_n) =1$, and set 
  $\eta_n = D^{\vec{v}_n}(\omega_n)$. Then for 
 $[\zeta]_{\F} \in H_c^{m-1}(M,\F) $  and each $n \geq 1$, we have
 \begin{equation}\label{eq-pairingetan}
G_{\F}(B)([\zeta]_{\F}) = \int_B \;  \eta_n \wedge \zeta ~ .
\end{equation}
Estimate the norms of the integrals in \eqref{eq-pairingetan}:  
\begin{eqnarray*}
\left| G_{\F}(B)([\zeta]_{\F}) \right| & = &  \lim_{n \to \infty} ~ \left| \int_B \;  \eta_n \wedge \zeta \right| \\
& \leq &  \lim_{n \to \infty} ~  \int_B \; \|  \eta_n \|_B \; \|\zeta \|_B  \; dvol \\
& \leq &  \lim_{n \to \infty} ~  (1/n) \cdot \int_B \; \|\zeta \|_B  \; dvol \\
& = & 0  ~ . 
\end{eqnarray*}
 As this holds for all  $[\zeta]_{\F} \in H_c^{m-1}(M,\F) $, the claim follows.
 \endproof
 
We note two important aspects of the   proof of Proposition~\ref{prop-vanishing}. First, the $n$-form $\eta_n \wedge \zeta$ in the integrand of \eqref{eq-pairingetan} depends only on the   restrictions $\eta_n | L$ for leaves $L$ of $\F$.
Thus, the pairing $\eta_n \wedge \zeta$   is well-defined when $\F$ is a $C^{\infty,1}$-foliation.  
Also, the convergence of the integral in \eqref{eq-pairingetan}  as $n \to \infty$ uses the Lebesgue dominated convergence theorem, and can be applied assuming only that   the form $\zeta  \in A_c^{m-1}(M, \F)$ is continuous. In particular, for a $C^{\infty,2}$-foliation the form $d\eta$ is continuous, so the calculation above applies to  multiples of this form as required for the proof of Proposition~\ref{prop-duality}.

 Proposition~\ref{prop-vanishing} gives an effective method for showing that the Godbillon-Vey class vanishes on a set $B \in \cB(\F)$, provided that one can construct a sequence of $1$-forms $\{\omega_n = \exp(f_n) \cdot \omega \mid n = 1, 2, \ldots\}$ on $M$ satisfying the hypotheses of the proposition. In hindsight, one can see that an analogous   estimate was  used in the   previous works 
 \cite{DuminySergiescu1981,Herman1977,  MMT1981, MoritaTsuboi1980, Tsuchiya1982,Wallet1976} 
to show that $GV(\F) = 0$ for $C^2$-foliations of codimension one, for foliations with various types of dynamical properties.
 
 For a  $C^2$-foliation $\F$,   Sacksteder's Theorem \cite{Sacksteder1965} implies that if $\F$ has no resilient leaf, then there are no exceptional minimal sets for $\F$.  Hence, by the Poincar\'e-Bendixson theory, all leaves of $\F$ either lie at finite level, or lie  in ``arbitrarily thin'' open subsets  $U \in \cB(\F)$. 
 In his works  \cite{Duminy1982a,Duminy1982b},    Duminy used a  result analogous to  Proposition~\ref{prop-vanishing} to show that $G_{\F}(B) = 0$, where $B$ is a union of leaves at finite level. 
  Thus, for a $C^2$-foliation with no resilient leaves, the Godbillon measure vanishes on the union of the leaves of finite level, and also vanishes on any Borel set in their complement.  Thus,  $GV(\F) = 0$ for a $C^2$-foliation of codimension-one with no resilient leaves. See    \cite{CandelConlon2000, CantwellConlon1984} for a published discussion of this proof.

In the next two sections, we follow a different, more direct approach   to obtain this conclusion. From  the assumption $G_{\F} \ne 0$, we conclude  that the holonomy pseudogroup of a $C^{\infty,1}$-foliation $\F$ must contain resilient orbits. Thus for  a $C^2$-foliation $\F$ with $GV(\F) \ne 0$, we have that  $G_{\F} \ne 0$  and hence  $\F$ must contain resilient leaves.

\section{Asymptotically expansive holonomy}\label{sec-pesin}

In this section,  we study the dynamical properties of $C^1$-pseudogroups acting on a $1$-dimensional space.
 The main example is when    there is given a codimension-one foliation $\F$ on a compact manifold $M$, with a regular $C^{\infty,1}$-foliation atlas with associated  transversal $\cX$, and $\GF$ is  its  holonomy pseudogroup. Then   $\GF$    is generated by a finite collection of local $C^1$-diffeomorphisms defined on open subsets of $\cX$. 
  Recall that the charts in the foliation atlas are assumed to be  transversally oriented, so for each plaque chain $\cP$,    the derivative    $\hh'_{\cP}(x)  > 0$ for all $x \in \cD_{\cP}$ in its domain. 

\subsection{The transverse expansion exponent function} 
We first introduce the notion of   \emph{asymptotically expansive}     holonomy for a  leaf of $\F$,  and the  associated set $\E$ of leaves with this property. 
  The main result of this section is  that  the Godbillon measure $G_{\F}$ is supported on  $\E$. That is, for any $B \in \cB(\F)$, we have 
 $G_{\F}(B) = G_{\F}(B \cap \E)$. 
Hence,   $G_{\F} \ne 0$ implies  the set $\E$ must have positive Lebesgue measure by Theorem~\ref{thm-duminy}.

For all $x \in \cX$, set $\mu_0(x) = 1$, and  and each integer $n \geq 1$, define the \emph{maximal n-expansion}
\begin{equation} \label{eq-est1}
\mu_n(x) =  \sup \, \{ \,  \hh'_{\cP}(x)  \mid x \in \cD_{\cP} ~ \& ~ \|\cP\| \leq n\} ~ .
\end{equation}
The function $x \mapsto \mu_n(x)$ is the maximum of a finite set of continuous functions, so is a Borel function on $\cX$, and  $\mu_n(x) \geq 1$ as 
 the identity transformation is the holonomy for a plaque chain of   length 1.

\begin{lemma}\label{lem-est}
Let $x \in \cX$,  and let  $\cQ = \{\cP_{\alpha}(x), \cP_{\beta}(y)\}$ be a plaque chain of length 1. For the holonomy map $\hh_{\cQ}$ of this length-one plaque-chain, we have   $\hh_{\cQ}(x) = y$.   Then  for all $n > 0$, 
\begin{equation} \label{eq-est3}
\mu_{n-1}(x)  \leq \mu_n(y)  \cdot \hh'_{\cQ}(x) \leq \mu_{n+1}(x) ~ .
\end{equation}
\end{lemma}
\proof
Let $\cP$ be a plaque chain at $y$ with $\|\cP\| \leq n$, then $\cP \circ \cQ$ is a plaque chain at $x$ with $\|\cP \circ \cQ\| \leq n+1$, so 
$$\hh_{\cP}'(y) \cdot \hh'_{\cQ}(x) = \hh_{\cP \circ \cQ}'(x)  \leq \mu_{n+1}(x) ~ .$$
As this is true for all plaque chains at $y$ with $\|\cP\| \leq n$, we obtain 
$\mu_n(y) \cdot \hh'_{\cQ}(x) \leq \mu_{n+1}(x)$.

Given a plaque chain $\cP$ at $x$ with $\|\cP\| \leq n-1$, the chain $\cR = \cP \circ \cQ^{-1}$ at $y$ has   $\|\cR\| \leq n$ and 
\begin{equation}\label{eq-est77}
  \hh_{\cP}'(x)  = \hh_{\cR}'(y) \cdot \hh'_{\cQ}(x)   \leq \mu_n(y)  \cdot \hh'_{\cQ}(x) ~ .
\end{equation}
As \eqref{eq-est77} holds  for all plaque chains at $x$ with $\|\cP\| \leq n-1$, we have 
$\mu_{n-1}(x)   \leq \mu_n(y)  \cdot \hh'_{\cQ}(x)$.
\endproof

\medskip

Define   $\lambda_n(x) =\ln \, (\mu_n(x))$, so that 
$\lambda_n(x)  = \sup \, \{  \ln(\hh'_{\cP}(x))   \mid x \in \cD_{\cP} ~ \& ~ \|\cP\| \leq n\}$.

Then the {\it transverse expansion exponent} at $x \in  \cX$ is   defined by 
\begin{equation} \label{eq-est2}
\lambda_*(x) =  \limsup_{n \to \infty}  \;  \frac{ \lambda_n(x) }{n} ~ .
\end{equation}

\begin{lemma}\label{lemma-inv}
The transverse expansion exponent function $\lambda_*$ is Borel measurable  on $\cX$, and  constant on the orbits of $\GF$.
\end{lemma}
\proof
For each $n \geq 1$, the function $\ds  \lambda_n(x)/n$ is Borel, so the supremum function in \eqref{eq-est2} is also Borel.

Let $x \in \cX$,  and let  $\cQ = \{\cP_{\alpha}(x), \cP_{\beta}(y)\}$ be a plaque chain, then   the estimate  (\ref{eq-est3}) implies that, 
\begin{equation}\label{eq-est2c}
\frac{ \ln (\mu_{n+1}(x)) }{n+1}  \geq  \frac{ \ln ( \mu_n(y)  \cdot \hh'_{\ba}(x)) }{n} \cdot \frac{n}{n+1} = \left\{\frac{ \ln ( \mu_n(y)) }{n}  + \frac{ \ln (  \hh'_{\ba}(x)) }{n}\right\}\cdot \frac{n}{n+1}
\end{equation}
so that 
\begin{equation}\label{eq-est2d}
\lambda_*(x) = \limsup_{n \to \infty}  \; \left\{ \frac{ \ln (\mu_{n+1}(x)) }{n+1} \right\} \geq   \limsup_{n \to \infty}  \left\{\frac{ \ln ( \mu_n(y)) }{n}\right\}= \lambda_*(y) ~ .
\end{equation}
The converse inequality follows similarly. 

Thus,   $\lambda_*(x) = \lambda_*(y)$  if there is a plaque chain $\cQ = \{\cP_{\alpha}(x), \cP_{\beta}(y)\}$. 
The pseudogroup   $\GF$ is generated by the holonomy defined by plaque chains of length $1$, so that for  each  point 
  $y \in \GF(x)$, there is a finite plaque chain  $\ds \cP = \{\cP_{\alpha_0}(x_0), \ldots, \cP_{\alpha_k}(x_k) \}$  with $x_0 = x$ and $x_k = y$. 
  Then  $\ds \lambda_*(x_{\ell}) = \lambda_*(x_{\ell + 1})$ for each $0 \leq \ell < k$, from which it follows that  $\lambda_*(x) = \lambda_*(y)$. 
 \endproof

\subsection{The expansion decomposition}
We use the conclusion of Lemma~\ref{lemma-inv} to lift the transverse expansion exponent function $\lambda_*$ from $\cX$ to $M$. 
For $\xi \in M$, let $L_{\xi}$ be the leaf containing $\xi$ and let $x = \pi_{\alpha}(\xi)$ where $\xi \in U_{\alpha}$. Then by Lemma~\ref{lemma-inv}, the value $\lambda_*(x)$ is independent of the choice of open set with $\xi \in U_{\alpha}$. By abuse of notation, we set   $\lambda_*(\xi) = \lambda_*(x)$, which is a well-defined function on $M$.
 Moreover, given a leaf $L$,  set 
$\ds \lambda_*(L) = \lambda_*(\xi)$ for some $\xi \in L$, which is then well-defined as well. 

\begin{defn}\label{defn-hc} Define the $\F$-saturated Borel subsets of $M$:
\begin{eqnarray*}
\E ~ & = & ~  \{ x \in  M \mid \lambda_*(x) > 0\}\\
\Ea ~ & = & ~  \{ x \in M \mid \lambda_*(x) > a\}, ~ {\rm for} ~ a  \geq  0\\
\SF ~ & = & ~  M - \E ~ .
\end{eqnarray*}
\end{defn}
 
A point $x \in \E$ is said to be {\em infinitesimally expansive}.  The set $\E$ is called the \emph{hyperbolic set} for $\F$, and  is the analog for codimension-one foliations of the hyperbolic set for diffeomorphisms in Pesin theory \cite{BarreiraPesin2013,Pesin1977}.   The set $\SF$  consists of the leaves  of $\F$ for which the transverse infinitesimal holonomy has    ``slow growth''.  Both    sets $\E$ and $\SF$ are fundamental for the study of the dynamics of the foliation $\F$.

Note that if there is an  holonomy map 
 $\hh_{\cP}$ with $ x \in \cD_{\cP}$,  $ \hh_{\cP}(x) = x$ and $ \hh'_{\cP}(x) = \lambda > 1$, then $x \in \E$. If $\cP$ is a plaque-chain of length $k$, then   $x \in \Ea$ for any $0 < a < \ln(\lambda)/k$.
 The plaque chain $\cP$ determines a    closed loop $\gamma_{\cP}$ based at $x$ in the leaf $L_x$, and the transverse holonomy along $\gamma_{\cP}$  is  linearly expanding in some open neighborhood of $x$.  Such  transversally  hyperbolic elements of the leaf holonomy   have a fundamental role in the study of foliation dynamics,  in particular in the  works by   Sacksteder \cite{Sacksteder1965}, by   Bonatti, Langevin and Moussu \cite{BonattiLangevinMoussu1992}, and the works \cite{Hurder1988,Hurder1991}. 
  However, given   $x \in \E$ there may not be a closed leafwise loop with infinitesimally expansive holonomy at $x$. What is always true is that   there is a sequence of holonomy elements whose length tends to infinity, what has infinitesimally expansive holonomy at $x$. We make this statement precise.

Consider a point  $x \in \cX \cap  \Ea$ for $a > 0$, and choose $\lambda$ with    $a < \lambda < \lambda_*(x)$.
Then for all $N > 0$, there exists   $n \geq N$ such that $ \lambda_n(x) \geq n \lambda$. By the definition of $\lambda_n(x)$, this means there exists a plaque chain $ \cP$ with length $ \|\cP\| \leq n$   starting  at $x$ such that   $\hh'_{\cP}(x)   \geq \exp \{n \lambda\}$. 
By the continuity of the derivative function on  $\cX$, there exists  $\e_n > 0$ such that  on the open interval $(x -\e_n, x+\e_n) \subset \cX$, 
$$ \hh'_{\cP}(y) \geq \exp \{  n \lambda/2\}   ~ {\rm for  ~  all} ~  x -\e_n \leq y \leq x+\e_n ~ .$$
By the Mean Value Theorem,   $\hh'_{\cP}$ is expanding on the interval $(x -\e_n, x+\e_n)$ by a factor at least $\exp \{  n \lambda/2\}$.
 Thus, the assumption $\lambda_*(x) > \lambda > 0$ and the definition in    \eqref{eq-est2} implies that we can choose a sequence of plaque 
chains $ \cP_{\ell}$ with lengths $ \|\cP_{\ell}\| = n_{\ell}$   starting  at $x$ such that  $n_{\ell}$   is strictly increasing, and so tends to infinity, and the corresponding holonomy maps satisfy
\begin{equation}\label{eq-holonomyexpansion}
 \hh'_{\cP_{\ell}}(y) \geq \exp \{  n_{\ell} \lambda/2\}   ~ {\rm for  ~  all} ~  x -\e_{n_{\ell}} \leq y \leq x+\e_{n_{\ell}} ~ .
\end{equation}
The constant   $\e_{n_{\ell}} > 0$ in \eqref{eq-holonomyexpansion} depends upon  $\ell$, $\lambda$ and $x$, and is exponentially decreasing  as $\ell \to \infty$.

It is  a strong condition to have a sequence of holonomy maps as in \eqref{eq-holonomyexpansion}   for elements of the holonomy pseudogroup at  points $x$, whose plaque lengths tend to infinity. This  is what gives the set $\E$ a fundamental role in the study of foliation dynamics, exactly in analog with the role of the Pesin set in smooth dynamics \cite{BarreiraPesin2013,Katok1980,Pesin1977,Ruelle1979}.
 The works \cite{Hurder2008,Hurder2010} give further study of the relation between the hyperbolic set $\E$ and the dynamics of the foliation.

In contrast, for the \emph{slow set}   $\SF$,    the dynamics of $\F$ on $\SF$ has ``less complexity'', as discussed in \cite{Hurder2010}. 
We next show that $G_{\F}(\SF) =0$, which is a measure of this lack of dynamical complexity.

Note that for an arbitrary saturated Borel set  $B \in \cB(\F)$, we have 
\begin{equation}\label{eq-decomposition}
 G_{\F}(B)  = G_{\F}(B \cap \E) + G_{\F}(B \cap \SF) 
\end{equation}
so  that $G_{\F} \ne 0$ and $G_{\F}(\SF) =0$ implies the set $\E$ must have positive Lebesgue measure.

\subsection{A vanishing criterion}
We use the criteria of  Proposition~\ref{prop-vanishing} to  show that $G_{\F}(\SF) =0$. That is,    we construct  a sequence of transversally measurable,  non-vanishing  transverse 1-forms   $\{\omega_n \mid n = 1, 2, \ldots\}$  on $M$ for which  $\| D^{\vec{v_n}}\omega_n\|_{\SF} < 1/n$.   
The construction of the forms $\{\omega_n \}$ follows the  method introduced in \cite{Hurder1986}. The first, and crucial step, is to construct an   $\e$--tempered cocycle (as given by \eqref{eq-tempered}) over the pseudogroup $\GF$  which is cohomologous to the  additive derivative cocycle,  using   a   procedure  adapted from \cite{HurderKatok1987}. This tempered cocycle is then used  to produce the sequence of defining 1-forms  $\omega_n$,   using the methods of  \cite{Bott1978}  and \cite{KT1973b, KT1975b}.
These are the used in the proof of the following result.
\begin{thm}\label{thm3} For any set $B \in \cB(\F)$, the Godbillon measure $G_{\F}(B) = G_{\F}(B \cap \E)$. \linebreak
Hence, if $\E$ has Lebesgue measure zero, then $G_{\F}(B) = 0$ for all $B \in \cB(\F)$. 
\end{thm}
\proof

By the above remarks, it suffices to show that $G_{\F}(\SF) =0$. 
Fix   $\e > 0$. 

For    $x \in  \cX \cap \SF$, by the definition of $\lambda_*(x) = 0$ \eqref{eq-est2},     there exists $N_{\e,x}$ such that 
$n \geq N_{\e,x}$ implies
$  \ln \{\mu_n(x)\} \leq n \e/2$, and hence the maximal $n$-expansion $\mu_n(x)  \leq \exp\{ n \e/2\}$. 

For $x \in \cX$ but  $x \not\in \SF$, set $g_{\e}(x) = 1$.  
For $x \in \cX \cap \SF$, set 
\begin{equation}\label{eq-lyapunov}
  g_{\e}(x) = \sum_{n=0}^{\infty}~ \exp\{ - n \e\} \cdot \mu_n(x) .
\end{equation}
  For $x$   in the slow set $\SF$, the sum in \eqref{eq-lyapunov}  converges    as the function $\exp\{ - n \e\} \cdot \mu_n(x)$ decays exponentially fast as $n \to \infty$. 
Note that while $g_{\e}(x) $ is finite for each $x \in \cX$, there need not be an upper bound for its values on $\cX \cap \SF$. 
Also,  $g_{\e}$ is  a Borel measurable function   defined on all of $\cX$.

The definition of the function $g_{\e}$    in \eqref{eq-lyapunov} is analogous to the definition of the Lyapunov metric in Pesin theory. Its role is to give a ``change of gauge''  with respect to which    the expansion rates of the dynamical system is ``normalized'' for the action of $\GF$ on $\cX$, as made precise by Lemma~\ref{lem-est2}  below.

 Let $\vec{X} = {\partial \over{\partial x}}$ be the unit-length, positively-oriented   vector field on $\mR$, let $dx$ denote the dual $1$-form on $\mR$.
 Recall that  for each $1 \leq \alpha \leq \nu$, we defined $\cX_{\alpha} \equiv (-1,1)$ so there is an inclusion $\iota_{\alpha} \colon \cX_{\alpha} \subset \mR$ which defines a coordinate function $x_{\alpha} \colon \cX_{\alpha} \to \mR$.   Then let   $dx_{\alpha} = \iota_{\alpha}^*(dx)$ denote the induced   $1$-form on  $\cX_{\alpha}$. There is a corresponding unit vector field $\vec{X}_{\alpha}$ on $\cX_{\alpha}$, which defines the vector field  $\vec{X}$ on $\cX$.  Then $dx_{\alpha}(\vec{X}) = 1$ on each $\cX_{\alpha}$.  
 
  For each $1 \leq \alpha \leq \nu$,  introduce the notation   $g_{\e}^{\alpha} = g_{\e} | \cX_{\alpha}$, and define   the  $1$-form  $dx_{\alpha}^{\e} = g_{\e}^{\alpha} \; dx_{\alpha}$ on $\cX_{\alpha}$.   

Let $x \in \cX$,  and let  $\cQ = \{\cP_{\alpha}(x), \cP_{\beta}(y)\}$ be a plaque chain of length $1$.  
Then  for the holonomy map $\hh_{\cQ}$ of this length-one plaque-chain, we have   $\hh_{\cQ}(x) = y$, and  
$\hh_{\cQ}^*(dx_{\beta}) = \hh_{\cQ}'  \cdot dx_{\alpha}$. Thus we have 
\begin{equation}\label{eq-coordchange}
\hh_{\cQ}^*(dx_{\beta}^{\e}) = \left( g_{\e}^{\beta} \circ  \hh_{\cQ}\right)  \cdot \hh_{\cQ}'  \cdot   dx_{\alpha}.
\end{equation}
The following result gives a key property of the function $g_{\e}$ which describes its behavior under a change of coordinates, for charts such that 
$\cP_{\alpha}(x) \cap \cP_{\beta}(y) \ne \emptyset$. Recall that $\SF = M - \E$ denotes the slow set, and let $S_{\alpha} = \pi_{\alpha}(\SF \cap U_{\alpha}) \subset \cX_{\alpha}$.

\begin{lemma}\label{lem-est2} 
For    $x \in  S_{\alpha}$  and  $\cQ = \{\cP_{\alpha}(x), \cP_{\beta}(y)\}$, 
\begin{equation}\label{eq-est4}
 \exp\{-\e\}  \cdot g_{\e}^{\alpha}(x)  \leq   g_{\e}^{\beta}(y) \cdot   \hh'_{\cQ}(x)  \leq  \exp\{\e\} \cdot g_{\e}^{\alpha}(x) ~ .
\end{equation}
\end{lemma}
\proof
Evaluate the expression  \eqref{eq-coordchange} on the vector field $\vec{X}$ and use the estimate (\ref{eq-est3}), noting that $\hh_{\cQ}(x) = y$,  to obtain, 
\begin{eqnarray} 
 g_{\e}^{\beta}(y) \cdot \hh_{\cQ}'(x) & = &  \left\{ \sum_{n=0}^{\infty} ~ \exp\{ - n \e\} \cdot \mu_n(y)\right\} \; \hh_{\cQ}'(x)    \nonumber \\
 & \leq &   \sum_{n=0}^{\infty} ~ \exp\{ - n \e\} \cdot \mu_{n+1}(x)  \nonumber \\
  & < &  \exp\{\e\}  \cdot \left\{ \sum_{n=1}^{\infty} ~ \exp\{ - n \e\} \cdot \mu_n(x)  + \mu_0(x) \right\} \nonumber \\
 & =   &  \exp\{\e\}  \cdot  g_{\e}^{\alpha}(x)  ~ . \label{eq-upper}
\end{eqnarray}

Similarly, we have 
\begin{eqnarray} 
 g_{\e}^{\beta}(y) \cdot \hh_{\cQ}'(x)   & = &  \left\{ \sum_{n=0}^{\infty} ~ \exp\{ - n \e\} \cdot \mu_n(y)\right\} \; \hh_{\cQ}'(x)    \nonumber \\
& \geq &   \sum_{n=1}^{\infty} ~ \exp\{ - n \e\} \cdot \mu_{n-1}(x)  + \mu_0(x) \cdot   \hh_{\cQ}'(x)  \nonumber \\
  & \geq &    \exp\{-\e\}  \cdot    ~ g_{\e}^{\alpha}(x) ~ .    \label{eq-lower}
\end{eqnarray}
This completes the proof of Lemma~\ref{lem-est2}.  
\endproof

 We next     use the coordinate $1$-form $dx_{\alpha}^{\e}$ to define a transversally measurable $1$-form $\omega_{\e}$ on $M$ which defines $\F$.
 The first step is to   define  local 1--forms $\omega^{\alpha}_{\e}$  on the coordinate charts  $U_{\alpha}$, then use a partition of unity to  obtain the $1$-form  $\omega_{\e}$ defined on all of $M$. Then, for  appropriate choices of $\e$ tending to $0$, we obtain  $1$-forms $\{\omega_n \}$ satisfying the hypotheses of    Proposition~\ref{prop-vanishing} on $\SF$.

 For each $1 \leq \alpha \leq \nu$, use the projection $\pi_{\alpha} \colon U_{\alpha} \to \cX_{\alpha}$ along plaques to pull-back  the form $dx_{\alpha}$ to the closed $1$-form $\omega_{\alpha} = \pi_{\alpha}^*(dx_{\alpha})$ on $U_{\alpha}$.  
Then define    $\ds \omega^{\alpha}_{\e}=  \pi_{\alpha}^*(dx_{\alpha}^{\e}) = \left( g_{\e}^{\alpha} \circ \pi_{\alpha}\right) \cdot \omega_{\alpha}$   which is a transversally measurable, leafwise closed  $1$-form on $U_{\alpha}$.

 Choose a partition of unity $\{\rho_{\alpha} \mid \alpha \in \cA\}$ 
subordinate  to the  cover $\{ U_{\alpha} \mid \alpha \in \cA\}$ of $M$ by foliation charts. Then for each $1 \leq \beta \leq \nu$, the $1$-form $\rho_{\beta} \cdot \omega^{\beta}_{\e}$ has support contained in $U_{\beta}$. 
Define the 1-form $\omega_{\e} = \sum \, \rho_{\beta} \cdot \omega^{\beta}_{\e}$ on $M$.  That is, for each $1 \leq \alpha \leq \nu$,  the restriction $\omega_{\e}| U_{\alpha}$  to the  chart $U_{\alpha}$ is given by 
\begin{equation}\label{eq-oneform}
 \omega_{\e}|{U_{\alpha}}   = \sum_{U_{\beta} \cap U_{\alpha} \not= \emptyset} \; \rho_{\beta}|_{U_{\alpha}}  \cdot   \omega^{\beta}_{\e} |_{U_{\alpha}} ~ .
\end{equation}

Recall that  $\vec{n}$ denotes the unit, positively-oriented vector field on $M$ orthogonal to $\F$, and let $\omega$ be the $1$-form on $M$ defining $\F$ with $\omega(\vec{n}) =1$.
Set $f_{\e} = \ln (\omega_{\e}(\vec{n}))$  so that $\omega_{\e} = \exp(f_{\e}) \cdot \omega$. Then for $\vec{v}_{\e} = \exp(- f_{\e}) \cdot \vec{n}$ we have 
$\omega_{\e}(\vec{v}_{\e}) =1$. 
Let 
\begin{equation}\label{eq-etae}
\eta_{\e}  = D^{\vec{v}_{\e}} \omega_{\e} =  \iota(\vec{v}_{\e}) d\omega_{\e} =\exp(- f_{\e}) \cdot   \iota(\vec{n}) d\omega_{\e}  ~ .
\end{equation}
 Recall that for the definition of the Godbillon measure,  we are only concerned   with  the restricted $1$-forms $\eta_{\e} | L$ for each leaf $L$  of $\F$. That is, the integrand in \eqref{eq-measure} depends only on the restricted class $\eta_{\e}|\F \in \Omega^1(\F)$. Recall from \eqref{eq-defDF} the definition of the leafwise differential   $\DF f_{\e} \equiv d(f_{\e}|\F) \in \Omega^1(\F)$. Then using  \eqref{eq-etae} we have the leafwise calculation in $\Omega^1(\F)$:  
\begin{eqnarray*}
\eta_{\e} |\F   & = &     \exp(- f_{\e}) | \F  \cdot   \left\{\iota(\vec{n}) d (\exp(f_{\e}) \cdot \omega)  \right\} |\F \\
  & = &    \exp(- f_{\e}) | \F  \cdot    \left\{  -\exp(f_{\e})| \F \cdot  \DF(f_{\e}) + \exp(f_{\e}|\F) \cdot \iota(\vec{n}) d \omega)|\F  \right\}  \\ 
  & = &  - \DF(f_{\e}) +  \left\{ \iota(\vec{n})  d \omega)  \right\} |\F \\  
&  = &   \eta |\F -  \DF(f_{\e}) ~ .
\end{eqnarray*}
Thus, the leafwise $1$-forms $\eta_{\e} |\F$ and $\eta|\F$   differ by the leafwise exact $1$-form $\DF(f_{\e})$. Then by the Leafwise Stokes' Theorem   \cite[Proposition~2.6]{HeitschHurder1984},  the Godbillon measure $G_{\F}(B)$ can be calculated using the   $1$-form  $\eta_{\e}|\F$ restricted to $B$.

We next estimate the norm  $\|\eta_{\e}\|$. Consider the $1$-forms $\ds \omega^{\beta}_{\e} |_{U_{\alpha}}$ appearing in the expression \eqref{eq-oneform}.
Fix $1 \leq \alpha \leq \nu$, then for $(\alpha, \beta)$ admissible, that is  such that $U_{\alpha} \cap U_{\beta} \not= \emptyset$,
let  $\cQ = \{\cP_{\alpha}(x), \cP_{\beta}(y)\}$ be a plaque chain with  holonomy map $\hh_{\cQ}$. Using the identity $ \pi_{\beta} =   \hh_{\cQ} \circ \pi_{\alpha}$ on 
  $U_{\alpha} \cap U_{\beta}$ and the identity \eqref{eq-coordchange}, then on $U_{\alpha} \cap U_{\beta}$  we have 
\begin{eqnarray*}
\omega^{\beta}_{\e} |_{U_{\alpha} \cap U_{\beta}} & = &  \pi_{\beta}^*(dx_{\beta}^{\e}) |_{U_{\alpha} \cap U_{\beta}} \\
  & = &  \pi_{\alpha}^* \circ \hh_{\cQ}^*(dx_{\beta}^{\e}) |_{U_{\alpha} \cap U_{\beta}}\\
  & = &  \pi_{\alpha}^*(g_{\e}^{\beta} \circ  \hh_{\cQ} \cdot \hh_{\cQ}'  \cdot   dx_{\alpha}) |_{U_{\alpha} \cap U_{\beta}} \\
  & = &  (g_{\e}^{\beta} \circ  \hh_{\cQ} \circ \pi_{\alpha})  \cdot (\hh_{\cQ}' \circ \pi_{\alpha}) \cdot \pi_{\alpha}^*(dx_{\alpha}) |_{U_{\alpha} \cap U_{\beta}} ~ .
\end{eqnarray*}

To simplify notation, set  $k_{\e,\alpha \beta}(x) = (g_{\e}^{\beta} \circ  \hh_{\cQ} \circ \pi_{\alpha})  \cdot (\hh_{\cQ}' \circ \pi_{\alpha})$. 
Note that $k_{\e,\alpha\alpha}= g_{\e}^{\alpha} \circ \pi_{\alpha}$.
Then in this notation, the estimate   \eqref{eq-est4} implies for $x \in U_{\alpha} \cap U_{\beta}$ that  
\begin{equation}\label{eq-tempered}
 \exp (-\e) \cdot k_{\e,\alpha \alpha}(x) \leq k_{\e,\alpha \beta}(x)   \leq   \exp (\e) \cdot k_{\e,\alpha \alpha}(x) ~.
 \end{equation}
Also, each function   $k_{\e,\alpha \beta}$ is constant along the plaques in $U_{\alpha} \cap U_{\beta}$, so that  its leafwise differential is zero; that is,  $\DF k_{\e,\alpha \beta} = 0$.
Recall that  $\omega_{\alpha} = \pi_{\alpha}^*(dx_{\alpha})$, so that $d\omega_{\alpha} = 0$, and  for $x \in U_{\alpha} \cap U_{\beta}$ we then have  
\begin{equation}\label{eq-cob}
\omega^{\beta}_{\e} |_x  = k_{\e,\alpha \beta}(x)  \cdot \omega_{\alpha}|_x ~ .
 \end{equation}
 
Then   for $x \in U_{\alpha}$ and using  the formulas \eqref{eq-oneform},  \eqref{eq-etae} and \eqref{eq-cob}, and letting $\vec{n}_x$ denote the value of the unit vector field $\vec{n}$ at $x$, we estimate $\left\| \eta_{\e} |_x \right\|$  as follows:
\begin{eqnarray} 
\left\| \eta_{\e} |_x \right\|    =    \left\| \exp(- f_{\e}(x)) \cdot   \{\iota(\vec{n}) d\omega_{\e}\}|_x \right\|    
&  = & \exp(- f_{\e}(x)) \cdot  \left\| \sum_{U_{\beta} \cap U_{\alpha} \not= \emptyset} \; \iota(\vec{n}_x) d \{\rho_{\beta}  \cdot   \omega^{\beta}_{\e} \}|_x  \right\| \nonumber \\   
&  = & \exp(- f_{\e}(x)) \cdot  \left\| \sum_{U_{\beta} \cap U_{\alpha} \not= \emptyset} \; \iota(\vec{n}_x) d \{\rho_{\beta}  \cdot   k_{\e,\alpha \beta} \cdot \omega_{\alpha}\} |_x \right\| \nonumber \\
&  = & \exp(- f_{\e}(x)) \cdot  \left\| \sum_{U_{\beta} \cap U_{\alpha} \not= \emptyset} \; \DF(\rho_{\beta})|_x  \cdot   k_{\e,\alpha \beta}(x) \cdot \omega_{\alpha}(\vec{n}_x)  \right\| ~ .\label{eq-almost1} 
\end{eqnarray} 

The the leafwise differential of the constant function is zero, so we have  the identity 
$$ 0 = \DF(1) = \DF(\sum \rho_{\beta})=  \sum \DF \rho_{\beta} .$$
We conclude that 
\begin{equation}\label{eq-vanishid}
0 = \sum_{U_{\beta} \cap U_{\alpha} \not= \emptyset} \DF(\rho_{\beta})|_x = \sum_{U_{\beta} \cap U_{\alpha} \not= \emptyset}  \DF(\rho_{\beta})|_x \cdot k_{\e,\alpha \alpha}(x) \cdot \omega_{\alpha}(\vec{n}_x) ~. 
\end{equation}
Then continuing from \eqref{eq-almost1},  and using the identities  \eqref{eq-vanishid} and \eqref{eq-tempered}, for $x \in U_{\alpha}$ we have:
\begin{eqnarray} 
\left\| \eta_{\e} |_x \right\|    
&  = &  \exp(- f_{\e}(x)) \cdot \left\| \sum_{U_{\beta} \cap U_{\alpha} \not= \emptyset} \; \DF(\rho_{\beta})|_x  \cdot   \{k_{\e,\alpha \beta}(x) - k_{\e,\alpha \alpha}(x)\} \cdot \omega_{\alpha}(\vec{n}_x) \right\| \nonumber \\
&  \leq &  \exp(- f_{\e}(x)) \cdot  \sum_{U_{\beta} \cap U_{\alpha} \not= \emptyset} \; \| \DF(\rho_{\beta})|_x \|  \cdot   | k_{\e,\alpha \beta}(x) - k_{\e,\alpha \alpha}(x)| \cdot | \omega_{\alpha}(\vec{n}_x)|  \nonumber \\
& \leq  &   \exp(- f_{\e}(x)) \cdot  \left\{ \sup_{x \in U_{\alpha}}  \| \DF(\rho_{\beta})|_x \| \cdot |\omega_{\alpha}(\vec{n}_x)| \right\} \cdot  \sum_{U_{\beta} \cap U_{\alpha} \not= \emptyset} \; | k_{\e,\alpha \beta}(x) - k_{\e,\alpha \alpha}(x)| \nonumber \\
& \leq  &   \exp(- f_{\e}(x)) \cdot  \left\{ \sup_{x \in U_{\alpha}}  \| \DF(\rho_{\beta})|_x \| \cdot |\omega_{\alpha}(\vec{n}_x)| \right\} \cdot  \sum_{U_{\beta} \cap U_{\alpha} \not= \emptyset} \; (\exp (\e) - 1) \cdot k_{\e,\alpha \alpha}(x)  ~ .  \label{eq-almost2} 
\end{eqnarray} 
It remains to estimate   $\exp(- f_{\e}(x))$   in \eqref{eq-almost2}.
Recall \eqref{eq-oneform} and using \eqref{eq-tempered} we have for $x \in U_{\alpha}$ that 
 \begin{eqnarray}
\exp(f_{\e}(x)) & = &      \sum_{U_{\beta} \cap U_{\alpha} \not= \emptyset} \; \rho_{\beta}(x)  \cdot   \omega^{\beta}_{\e}(\vec{n}_x)  \nonumber\\
& = &   \sum_{U_{\beta} \cap U_{\alpha} \not= \emptyset} \; \rho_{\beta}(x)  \cdot   k_{\e,\alpha \beta}(x) \cdot \omega_{\alpha}(\vec{n}_x) \nonumber\\
& \geq &   \sum_{U_{\beta} \cap U_{\alpha} \not= \emptyset} \; \rho_{\beta}(x)  \cdot    \exp (-\e) \cdot k_{\e,\alpha \alpha}(x) \cdot \omega_{\alpha}(\vec{n}_x) \nonumber \\
& = &   \exp (-\e) \cdot k_{\e,\alpha \alpha}(x) \cdot \omega_{\alpha}(\vec{n}_x)  ~ .\label{eq-almost3} 
\end{eqnarray}
Thus, we obtain the estimate
\begin{equation}\label{eq-upperest}
\exp(- f_{\e}(x)) \leq \exp (\e) \cdot (k_{\e,\alpha \alpha}(x)  \cdot \omega_{\alpha}(\vec{n}_x))^{-1} .
\end{equation}

Then combining \eqref{eq-almost2}  and \eqref{eq-upperest}, and noting that the number of indices $\beta$ for which  $U_{\beta} \cap U_{\alpha} \not= \emptyset$ is bounded by the cardinality $\nu$ of the covering,  we obtain
 \begin{equation}
\left\| \eta_{\e} |_x \right\|   \leq \left\{ \sup_{x \in U_{\alpha}}  \| \DF(\rho_{\beta})|_x \|   \right\} \cdot  \nu \cdot  \exp(\e) (\exp (\e) - 1)   
   \label{eq-final} 
\end{equation}

Note that the right hand side in \eqref{eq-final} tends to 0 as $\e \to 0$, so that 
for  each $n > 0$,  we can  choose $\e_n > 0$ such that $\|\eta_{\e_n}\| \leq 1/n$. Then  set $\omega_n = \omega_{\e_n}$, and the claim of the Theorem~\ref{thm3} follows.
\endproof

\section{Uniform hyperbolic expansion}\label{sec-uniform}

 In this section, we assume  that $\F$ is a $C^1$-foliation with non-empty hyperbolic set $\E$, and  show that  there exists a hyperbolic fixed-point for the holonomy pseudogroup $\GF$. The proof uses a pseudogroup version of the \emph{Pliss Lemma}, which is fundamental in the study of non-uniformly hyperbolic dynamics (see \cite{ABV2000} or \cite[Lemma~11.5]{BDV2005}, or the original article by Pliss \cite{Pliss1972}.) 
 
 The goal is to construct hyperbolic contractions in the holonomy pseudogroup. The length of the path defining the holonomy element is not important, but rather it is important to obtain uniform  estimates on the size of the domain of the hyperbolic element   thus obtained, estimates which are independent of the length of the path. This is a key technical point   for the application of the constructions of this section in the next Section~\ref{sec-hyperbolic}, where we construct sufficiently many contractions so that they result in the existence of a resilient orbit for the action of the holonomy pseudogroup.

 We note that the existence of a hyperbolic contraction can also be deduced using the foliation geodesic flow methods introduced in \cite{Hurder1991}, though that method does not yield estimates on the size of the domain of the hyperbolic element in the foliation pseudogroup.  

\subsection{Uniform hyperbolicity and the Pliss Lemma}
We fix a regular covering on $M$ as in Section~\ref{subsec-regular}, with transversals $\fX$ and $\tX$ as in \eqref{eq-transversalX}, and let   $\GF$ denote the resulting pseudogroup acting on the spaces $\cX$ and $\tcX$ as in \eqref{eq-transversalT}.
Recall that by   Lemma~\ref{lem-leb}, 
there exists $\e_0 > 0$ so that for every admissible pair $(\alpha, \beta)$ and $x \in \cX_{\ab}$ then   $[x-\e_0, x+ \e_0]  \subset \tcX_{\ab}$. 
Recall that the space $\cX_{\ab}$ was defined in \eqref{eq-chiab}, and $\tcX_{\ab}$ was defined in \eqref{eq-wtchiab}.

\begin{defn}\label{def-delta0}
Given $0 < \e_1 \leq \e_0$, a constant   $0 < \delta_0 \leq \e_1$ is said to be a \emph{logarithmic modulus of continuity} for $\GF$ with respect to $\e_1$,  if      for  $y,z \in \cX_{\ab}[\delta_0]$   with   $\dT(y,z) \leq \delta_0$,     then  
\begin{equation}\label{eq-uniform}
 \left| \log\{\wh_{\ba}'(y)\} - \log\{\wh_{\ba}'(z)\} \right| \leq \e_1 ~ .
\end{equation}
\end{defn}

\begin{lemma}\label{lem-delta0}
Given $0 < \e_1 \leq \e_0$, there exists a constant   $0 < \delta_0 \leq \e_1$ which is a logarithmic modulus of continuity for $\GF$ with respect to $\e_1$.
\end{lemma}
\proof
By the choice of $0 < \e_1 \leq \e_0$, for each admissible pair $\{\alpha, \beta\}$,   the logarithmic derivative    $\log \{ \wh_{\ba}'(y)\}$ is continuous on the compact subset $\cX_{\ab}[\e_1] \subset \tcX_{\ab}$.  Thus, there exists $\delta_0(\alpha, \beta) > 0$ such that \eqref{eq-uniform} holds for this choice of $\{\alpha, \beta\}$. Define 
$\ds \delta_0 = \min \{\delta_0(\alpha, \beta) \mid \{\alpha, \beta\} ~ {\rm admissible}\}$. As the number of admissible pairs is finite, we have $\delta_0 > 0$. 
\endproof

 The next result shows  that if $\E$ is non-empty, then  there are words in $\GF$ of arbitrarily long length,  along which the  holonomy is ``uniformly expansive''. That is, there exists a constant $\lambda_* > 0$ such that for such a word   $\hh_n$ defined by a plaque chain $\cP$ of length $n$, then $\hh_n'(y) \geq \exp\{n \lambda_*)$ for all $\ds y \in \cD_{\cP}$.  
The proof is technical, but also notable  as it develops  a version for pseudogroup actions of   the    Pliss Lemma, which is used in the study of the dynamics of   partially hyperbolic diffeomorphisms, as for example  in \cite{BDV2005,Mane1987,Pliss1972}.    
 
Note that Definition~\ref{defn-hc} implies that the set $\E$ is an increasing union  of the sets $\Ea$   for   $a > 0$, and 
thus given   $\xi \in \E$,   there  exist $a> 0$ such that  $\xi \in \Ea$. 
 
We introduce a convenient notation  for working with the set $\Ea$. For each $1 \leq \alpha \leq \nu$, let 
\begin{eqnarray*}
\Ea \cap \cX_{\alpha} & = &  \pi_{\alpha}(\Ea \cap U_{\alpha}) \subset   \cX_{\alpha} \\
\Ea \cap \cX & = &  (\Ea \cap \cX_{1}) \cup \cdots \cup  (\Ea \cap \cX_{\nu}) .
\end{eqnarray*}
Recall that the transversals $\fX_{\alpha}$   and their images   $\cX_{\alpha}$ in the coordinates $U_{\alpha}$ were defined in \eqref{eq-transversals}.

\begin{figure}[htbp]
\begin{center}
\includegraphics[width=0.6\textwidth]{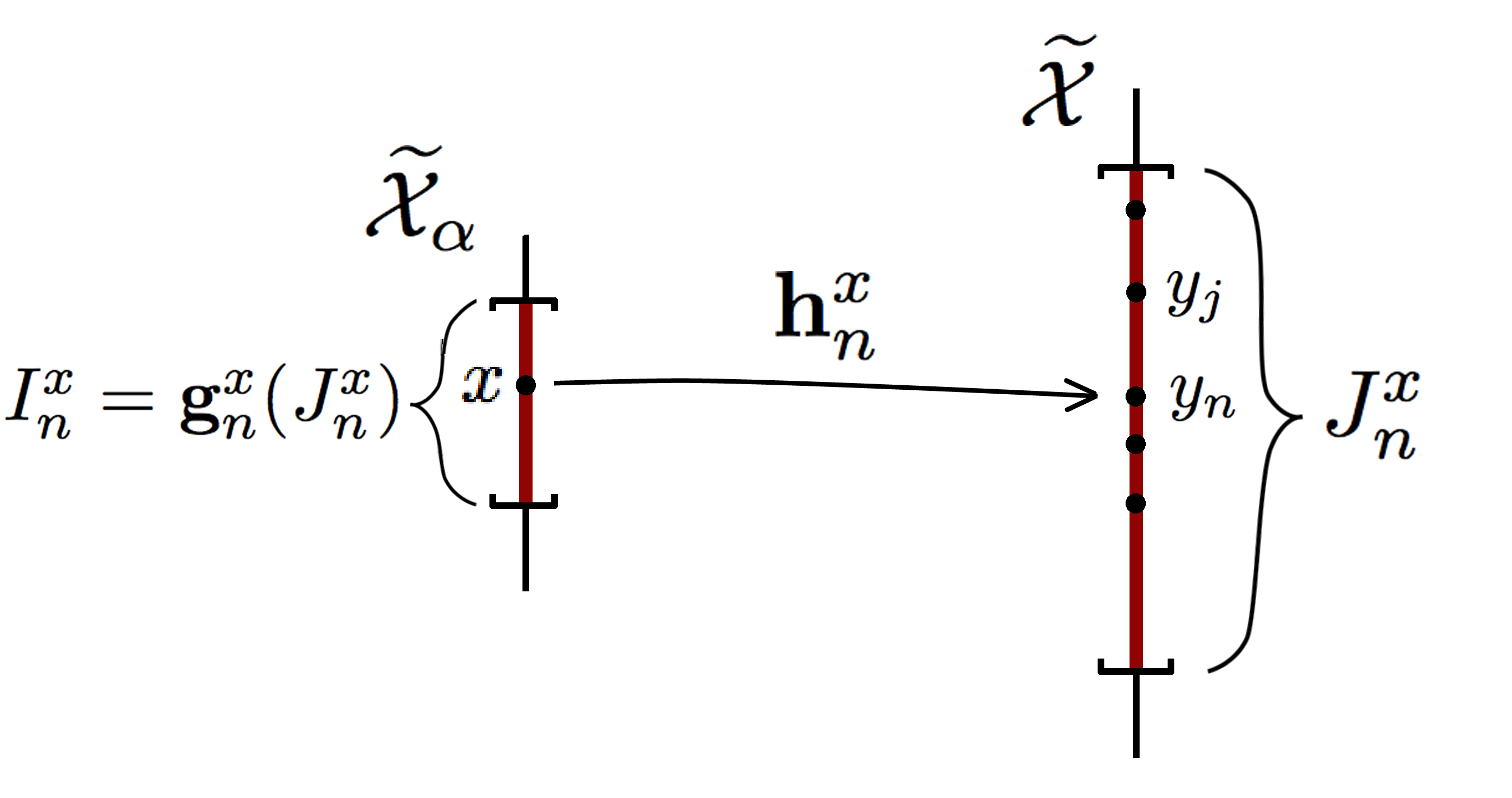}
\caption{Expanding holonomy map $\hh_n^x$}
\label{Figure4}
\end{center}
\end{figure}

\begin{prop}\label{prop-localhol}
Let $x   \in \Ea \cap \cX$ for $a > 0$,   let  $0 < \e_1 < \min\{\e_0, a/100\}$, and let $\delta_0$ be the logarithmic modulus of continuityÊ   for $\GF$ with respect to $\e_1$, as chosen   in Lemma~\ref{lem-delta0}. 

Then for each integer $n > 0$, there exist a point $y_n \in \GF(x)$, a closed interval $I_n^x   \subset \tcX_{\alpha}$  containing $x$ in its interior, and a  holonomy map  $\hh_n^x \colon I_n^x \to J_n^x$  such that for $y_n = \hh_n^x(x)$,  
  $\ds J_n^x = [y_n -  \delta_0/2, y_n+\delta_0/2] \subset \tcX$  and $I_n^x = (\hh_n^x)^{-1}(J_n^x)$,   
we have
\begin{equation}\label{eq-hypcontraction}
 (\hh_n^x)'(z) > \exp\{n a/2\}~ {\rm   for ~  all}  ~   z \in I_n^x ~ .
\end{equation}
It follows that $\ds |I_n^x| <  \delta_0 \, \exp\{-n a/2\}$. This is illustrated in  Figure~\ref{Figure4}.
\end{prop}

\proof
Fix a choice of   $0 < \e_1 < \min\{\e_0, a/100\}$, and then choose a logarithmic modulus of continuityÊ $\delta_0 > 0$ as in Lemma~\ref{lem-delta0}.

The set $\cX_{\ab}[\e_1]$, as defined in \eqref{eq-abdeltanbhd} for $\delta = \e_1$, is compact, so there exists  $C_0 > 0$ so that for all $(\alpha, \beta)$ admissible and  $y \in \cX_{\ab}[\e_1]$,  we have 
$\ds 1/C_0 \leq \wh_{\ba}'(y) \leq C_0$.

From the definition of $\lambda_*(x)$ as a $\limsup$ in \eqref{eq-est2}, the assumption that  $\lambda_*(x) > a$ implies that  for each   integer $n > 0$,
we can choose a   plaque chain of length $\ell_n \geq n$, given by  $\ds \cP_n = \{\cP_{\alpha_0}(z_0), \ldots, \cP_{\alpha_{\ell_n}}(z_{\ell_n}) \}$  with      $z_0 = x$,     such that    $\ds \log\{\hh_{\cP_n}'(z_0)\} >  \ell_n \cdot a$.  
Fix $n$ and the choice of the  plaque chain $\cP_n$ as above.

For each $1 \leq j \leq \ell_n$ let $ \hh_{\alpha_{j},\alpha_{j-1}}$ be the holonomy transformation defined by   $\{\cP_{\alpha_{j-1}},\cP_{\alpha_j}\}$, and so $ \hh_{\alpha_{j},\alpha_{j-1}}^{-1} =  \hh_{\alpha_{j-1},\alpha_{j}}$. 
Introduce the notation  $\whh_{0} = Id$, and   for $1 \leq j \leq \ell_n$ let
\begin{equation}\label{eq-hk}
\whh_{j}   =   \hh_{\alpha_{j},\alpha_{j-1}} \circ \cdots \circ \hh_{\alpha_{1},\alpha_{0}} 
\end{equation}
denote the partial composition of generators.  Note  that $z_j = \whh_{j}(z_0)$ and $z_0 = x$, and   that we have the   relations 
$\ds \whh_{j+1} = \hh_{\alpha_{j+1},\alpha_{j}} \circ \whh_{j}$ and   $z_{j+1} =   \whh_{j}(z_{j})$  for $0 \leq j < \ell_n$.
For  each $1 \leq j \leq \ell_n$, set 
\begin{equation}\label{eq-lambdaj}
\lambda_j = \log\{\whh'_{\alpha_{j-1},\alpha_{j}}(z_j)\} = - \log\{\whh'_{\alpha_{j},\alpha_{j-1}}(z_{j-1})\}  ~ .
\end{equation}
In particular,  $\ds \log\{\whh'_{\ell_n}(x)\} = -(\lambda_1 + \cdots + \lambda_{\ell_n})$. 
Note that if $\lambda_j < 0$ then the map $\whh_{\alpha_{j-1},\alpha_{j}}$ is an infinitesimal  contraction at $z_j$, and $\whh_{\alpha_{j},\alpha_{j-1}}$ is an infinitesimal  expansion at $z_{j-1}$.

The following algebraic definition and     lemma provide 
 the key to the analysis of the hyperbolic expansion properties of the partial compositions of the   maps $\whh_{j}$.

\begin{defn}\label{def-regular}
Let $\{\lambda_1, \ldots ,\lambda_m\}$ be given, and $\vartheta  > 0$.   An index $1 \leq j \leq m$ is said to be  \emph{$\vartheta$-regular} if   the following sequence of partial sum estimates hold: 
\begin{eqnarray}
\lambda_{j} + \;\; \vartheta & < &  0 \nonumber\\
\lambda_{j-1} + \lambda_{j} + 2 \vartheta & < & 0\nonumber\\
 & \vdots &  \label{eqs-regular}\\
\lambda_{1} + \cdots + \lambda_{j} + j \vartheta & < & 0 ~ .\nonumber
\end{eqnarray}
\end{defn}
 Condition \eqref{eqs-regular}   is a weaker hypothesis than   assuming the uniform estimates $\lambda_i <-\vartheta$ for all $1 \leq i \leq j$, but  is sufficient    for our purposes. 
The next result shows that $\vartheta$-regular indices always exist.

\begin{lemma}\label{lem-regular}
Assume   there are given  real numbers $\{\lambda_{1}, \ldots , \lambda_{m}\}$ such that 
\begin{equation}\label{eq-regular}
\lambda_{1} + \cdots + \lambda_{m}  \leq  -a \, m  ~ .
\end{equation}
 Then   for any  $0 < \e_1 < a$,  there exists an $\e_1$-regular index $q_m$, for some $1 \leq q_m \leq m$, which satisfies
\begin{equation}\label{eq-maxlem}
 \lambda_{1} + \cdots + \lambda_{q_m} ~ \leq ~ (-a + \e_1) \, m ~ .
\end{equation}
\end{lemma}
\proof
 
 The existence of the index $q_m$ satisfying this property is shown by contradiction. We introduce the concept of an \emph{$\e_1$-irregular} index, for which   the \emph{$\e_1$-regular}  condition fails, and show by contradiction that not all indices can be $\e_1$-irregular.

We say that an index $k \leq m$ is  \emph{$\e_1$-irregular}  if  
\begin{equation}\label{eq-irregular}
\lambda_{k} + \cdots + \lambda_{m} + (m - k + 1) \e_1 \geq 0 ~ .
\end{equation}
If there is no irregular  index, then observe that  $q_m = m$ is an $\e_1$-regular index.
 Otherwise, suppose that there exists  some   index $k$ which is   $\e_1$-irregular.
 The   inequality \eqref{eq-regular} states that  
   the   index $k= 1$ is not $\e_1$-irregular.  
Let $j_m \leq m$ be  the least     $\e_1$-irregular index, so that
\begin{equation}\label{eq-irregular2}
\lambda_{j_m} + \cdots + \lambda_{m} + (m - j_n + 1) \e_1 \geq 0
\end{equation}
By \eqref{eq-regular}, $j_m =1$ is  is not $\e_1$-irregular, so we have $2 \leq j_m \leq m$.

Set $q_m = j_m -1$, then we claim that   $q_m $ is an $\e_1$-regular index. 
If not, then at least one of the inequalities in \eqref{eqs-regular} must fail to hold. That is,  there is some $i \leq q_m$ with 
\begin{equation}\label{eq-patch}
\lambda_{i} + \cdots  + \lambda_{q_m} + (q_m-i +1) \e_1 \geq 0 ~ .
\end{equation}
Add the inequalities   (\ref{eq-irregular2}) and   (\ref{eq-patch}), and noting that $q_m = j_m -1$,  we obtain that  $i$ is also an $\e_1$-irregular index. As $i < j_m$, this is contrary to the choice of $j_m$.   Hence,   $q_m $ is an $\e_1$-regular index.

It remains to show that the estimate \eqref{eq-maxlem} holds.  As  $j_m = q_m +1$ is irregular,   subtract  (\ref{eq-irregular}) for $k = j_m$  from (\ref{eq-regular}) to obtain  
$$
 \lambda_{1} + \cdots + \lambda_{q_m}  \leq -a m + (m - q_m) \e_1 
 \leq (-a + \e_1) \, m   
$$
as claimed. 
\endproof

We return to considering the maps $\whh_j$ defined by \eqref{eq-hk}, and the exponents $\lambda_j$ defined by \eqref{eq-lambdaj}. The following result then follows directly from Lemma~\ref{lem-regular} and the definitions.

 \begin{cor}\label{cor-regular2}
 Assume that there is given   $a > 0$ with  $x \in \Ea \cap \cX$, a choice of integer $n > 0$, and plaque-chain 
$\ds \cP_n = \{\cP_{\alpha_0}(z_0), \ldots, \cP_{\alpha_{\ell_n}}(z_{\ell_n}) \}$ with $\ell_n \geq n$, 
 such that    $\ds \log\{\hh_{\cP_n}'(z_0)\} \geq  \ell_n \cdot a$  and  $z_0 = x$. 
 Given     $0 < \e_1 < a$,   by Lemma~\ref{lem-regular}  there exists an $\e_1$-regular index $q_n$, for some $1 \leq q_n \leq \ell_n$ chosen as in Lemma~\ref{lem-regular}, such that for the map $\whh_{q_n}$ defined by \eqref{eq-hk}, 
\begin{equation}\label{eq-contracting}
\whh_{q_n}'(x)  \geq (a - \e_1) \, \ell_n \geq (a - \e_1) \, n~ .
\end{equation}
 \end{cor}

The  estimate (\ref{eq-contracting})   can be interpreted as stating that ``most'' of   the  infinitesimal expansion of the map 
$\whh_{\ell_n}$ at $z_0$ is achieved by  the action of  the partial composition $\whh_{q_n}$.

Recall that we have  a fixed choice of   $0 < \e_1 < \min\{\e_0, a/100\}$, as given in the statement of Proposition~\ref{prop-localhol}, and  $\delta_0 > 0$ is chosen so that   the uniform continuity estimate  (\ref{eq-uniform}) in Lemma~\ref{lem-delta0} is satisfied. 

Then let $1 \leq q_n \leq \ell_n$ be the $\e_1$-regular index defined in Lemma~\ref{lem-regular} which satisfies \eqref{eq-maxlem}. 
We  next use  the $\e_1$-regular condition to obtain uniform estimates on the domains for which the inverses $\whh_{j}^{-1}$ are contracting, for $1 \leq j \leq q_n$.

Recall that $\ds \wh_{\alpha,\beta}$ denotes the continuous extension of the map $\ds \hh_{\alpha,\beta}$ to the   domain $\tcX_{\ab}$. 
Introduce       extensions  $\hh_n^x$ of $\whh_{q_n}$ and $\hg_n^x$ of its inverse $\whh_{q_n}^{-1}$, which  are defined by 
\begin{eqnarray}
\hh_n^x  & = & \wh_{\alpha_{q_n},\alpha_{q_n -1}} \circ \cdots \circ \wh_{\alpha_{1},\alpha_{0}}  \label{defn-hn}\\
\hg_n^x   & = &   \wh_{\alpha_{0},\alpha_{1}} \circ \cdots \circ \wh_{\alpha_{q_n-1},\alpha_{q_n}} ~ . \label{defn-gn}
\end{eqnarray}

Set  $y_n = \hh_n^x(x) = z_{\ell_n}$, then  by the estimate (\ref{eq-contracting}) we have 
\begin{equation}\label{eq-keyest}
 \log\{ (\hg_n^x)'(y_n)\}  ~ = ~  \lambda_{1} + \cdots + \lambda_{q_n} ~  \leq ~   (-a+\e_1) \; \ell_n < 0~.
 \end{equation}

We next show that  $\hg_n^x $ is uniformly contracting on an interval   with uniform length about $y_n$.

\begin{lemma}\label{lem-ie-domain}
Set $\delta_0' = \delta_0/8$.
Then the   interval $J_n^x = [y_n - 4\delta_0', y_n + 4\delta_0']$  is in the domain of $\hg_n^x$,   and for all $y \in J_n^x$, 
\begin{equation}\label{eq-contract1}
 \exp\{(-a - 2 \e_1)\, \ell_n\} \leq (\hg_n^x)'(y) \leq    \exp\{(-a + 2\e_1) \, \ell_n\} ~. 
\end{equation}
Hence, for $I_n^x = \hg_n^x(J_n^x)$,
\begin{equation}\label{eq-contract2}
 |I_n^x| \leq   \delta_0 \exp\{(-a + 2\e_1)\; \ell_n \} < \exp\{(-a/2)\; \ell_n \} ~ . 
\end{equation}
\end{lemma}
\proof
By the   choice of $\delta_0'$, the uniform continuity estimate  (\ref{eq-uniform})    implies that for all  $y  \in J_n^x$ 
$$ \left| \log\{ \wh'_{\alpha_{q_n -1},\alpha_{q_n}}(y)\} - \log\{\wh'_{\alpha_{q_n - 1},\alpha_{q_n}}(y_n)\} \right| \leq \e_1 ~ .
$$
Thus, by the definition of $\lambda_{q_n}$ we have that, for all $y \in J_n^x$,
$$ \exp\{\lambda_{q_n} -\e_1 \} \leq 
\wh'_{\alpha_{q_n - 1},\alpha_{q_n}}(y) \leq   
\exp\{\lambda_{q_n} + \e_1 \} ~ .$$ 
The assumption that  $q_n$ is $\e_1$-regular implies 
$\ds   \lambda_{q_n}+\e_1 < 0$, hence  
$\exp\{\lambda_{q_n}+ \e_1\} < 1$. 
Thus,  for all  $y \in J_{n}^x$ we have 
\begin{equation}\label{eq-ind1}
    \dT(\wh_{\alpha_{q_n -1},\alpha_{q_n}}(y_n), \wh_{\alpha_{q_n -1},\alpha_{q_n}}(y))    \leq 4\delta_0' \exp\{ \lambda_{q_n} + \e_1\} < 4 \delta_0' ~ .
 \end{equation}

Now proceed by downward  induction. For $0 <   j \leq q_n$ set 
$$\hg_{n,j}^x   =     \wh_{\alpha_{j-1},\alpha_{j}} \circ \cdots \circ \wh_{\alpha_{q_n -1},\alpha_{q_n}} ~, \;\;  J_{n,j}^x =  \hg_{n,j}^x(J_{n}^x) ~, \;\;  y_{n,j} = \hg_{n,j}^x(y_n) = z_{j-1}  ~ .$$

 Assume that for $1 <  j \leq q_n$,  we are given that for all 
$y \in J_{n,j}^x$ the   estimates 
\begin{equation}  \label{eq-induction1}
 \exp \{ \lambda_{j}  + \cdots + \lambda_{q_n} - (q_n -j +1) \, \e_1 \}  \leq   (\hg_{n,j}^x)'(y)   \leq   \exp \{ \lambda_{j}  + \cdots + \lambda_{q_n} + (q_n -j +1) \, \e_1 \} ~,  
\end{equation}
\begin{equation}\label{eq-induction2}
  \dT(y, y_{n,j})   \leq     4\delta_0' ~ .  
\end{equation}

The choice of $\delta_0$ and the   hypothesis (\ref{eq-induction2})  imply that for $y  \in J_{n,j}^x$,
$$ 
\left| \log\{ \wh'_{\alpha_{j -2},\alpha_{j-1}}(y)\} - \log\{ \wh'_{\alpha_{j -2},\alpha_{j-1}}(y_{n,j})\} \right| \leq \e_1 ~ .
$$

Recall that $z_{j-1} = y_{n,j}$, and  that $\lambda_{j-1} = \log\{ \wh'_{\alpha_{j -2},\alpha_{j-1}}(y_{n,j})\}$ by   \eqref{eq-lambdaj},    so   for all $y \in J_{n,j}^x$ we have for the inverse map $\wh_{\alpha_{j -2},\alpha_{j-1}} = \wh_{\alpha_{j -1},\alpha_{j-2}}^{-1}$ that 
\begin{equation}\label{eq-indstep}
 \exp\{\lambda_{j-1} -\e_1 \}   \leq 
\wh'_{\alpha_{j -2},\alpha_{j-1}}(y) \leq   
\exp\{\lambda_{j-1} + \e_1 \}   ~ .
\end{equation}

Then by the chain rule, the estimates (\ref{eq-indstep})    and the inductive hypothesis (\ref{eq-induction1})  yield the estimates
\begin{equation}
 \exp \{ \lambda_{j-1}  + \cdots + \lambda_{q_n} - (q_n -j +2) \, \e_1 \}  \leq   (\hg_{n,j-1}^x)'(y)   \leq   \exp \{ \lambda_{j-1}  + \cdots + \lambda_{q_n} + (q_n -j +2) \, \e_1 \}.  \label{eq-induction3}
\end{equation}

Now the assumption that  $q_n$ is $\e_1$-regular implies 
$\ds \lambda_{j-1}  + \cdots + \lambda_{q_n} + (q_n -j +2) \, \e_1 < 0$
 hence  
$ \exp\{\lambda_{j-1}  + \cdots + \lambda_{q_n} + (q_n -j +2) \, \e_1\} < 1$.

By the Mean Value Theorem, this yields  the distance bound        $\ds \dT(y_{n,j-1} , y)   \leq     4\delta_0'$, which is the hypothesis (\ref{eq-induction2})  for $j-1$. This completes the inductive step.  Thus, we may take $j=1$ in inequality (\ref{eq-induction1}) and   combined with the inequality   (\ref{eq-maxlem}),   for all 
$y \in J_{n}^x$ we have that 
\begin{equation}\label{eq-conjugatehyp}
(\hg_{n}^x)'(y)    \leq     \exp \{ \lambda_{1}  + \cdots + \lambda_{q_n} +  q_n  \, \e_1 \}    \leq    \exp \{ -a \, \ell_n  + (\ell_n + q_n ) \, \e_1 \}   
  \leq     \exp \{ (-a + 2\e_1) \, \ell_n    \}  ~ .
\end{equation}
Set  $I_n^x = \hg_n^x(J_n^x)$, then the estimate (\ref{eq-contract2})   follows by the Mean Value Theorem. 
\endproof

Since $ a - 2 \e_1 > a/2$ and $\ell_n \geq n$, this completes the proof of  Proposition~\ref{prop-localhol}. 
\endproof

 \subsection{Hyperbolic fixed-points}
We   show the existence of hyperbolic fixed-points for $\GF$ contained in the closure of   $\E \cap \cX$ in  $\tcX$, with uniform estimates  on the lengths of their  domains of contraction.

\begin{prop}\label{prop-fp}
Let $x   \in \Ea \cap \cX$ for $a > 0$,  
let  $0 < \e_1 < \min\{\e_0, a/100\}$, and let $\delta_0$ be chosen as in Lemma~\ref{lem-delta0}, and set $\delta_0' = \delta_0/8$.
Given $0 < \delta_1 < \delta_0'$ and $0 < \mu < 1$, then     there exists holonomy maps $\phi_1, \psi_1 \in \GF$, points $u_1, v_1 \in \overline{\cX}$ such that $\dT(x,v_1) < \delta_1$, such that we have:
 \begin{enumerate}
\item $\Phi_1 = \phi_1 \circ \psi_1$ has fixed point $\Phi_1(u_1) = u_1$;  
\item  $\cJ_1 \equiv [u_1 - \delta_0', u_1+ \delta_0']$ is contained in the domain of $\Phi_1$;  
\item $\Phi'(y) < \mu$ for all $y \in \cJ_1$; 
\item $\Psi_1 = \psi_1 \circ \phi_1$ has fixed point  $\Psi_1(v_1) = v_1$; 
\item  $\cK_1 \equiv \psi_1(\cJ_1) \subset (x - \delta_1, x+ \delta_1)$.
\end{enumerate}
\end{prop}
\proof
 The idea of the proof is to consider a sequence of maps as given by Proposition~\ref{prop-localhol}, for $n \geq 1$, and consider a subsequence of these for which the sequence of points $\ds \{y_n = \hh_n^x(x) = z_{\ell_n} \mid n \geq 1\}$ cluster at a limit point. We then use the estimates \eqref{eq-contract2} on the sizes of the domains 
to show that the appropriate compositions of these maps are defined, and have a hyperbolic fixed point. The details of this argument follow.

 Set $\delta_* = \min\{1, \delta_0'/4,\delta_1/4\}$.   Then by Proposition~\ref{prop-localhol}, for each integer $n > 0$,  we can choose a map  $\hh_n^x \colon I_n^x \to J_n^x$  as in \eqref{defn-hn}, which satisfies   condition  (\ref{eq-hypcontraction}). Label the resulting sequence of points   $y_n = \hh_n^x(x)  \in \cX$, and   the    inverse maps  $\hg_n^x = (\hh_n^x)^{-1}$.  Let $p_n$ denote the   length of the plaque chain defining $\hh_n^x$, then $p_n$ equals  the   $\e_1$-regular index $1 \leq q_n \leq \ell_n$ chosen as in the proof of Corollary~\ref{cor-regular2}.

 Recall that $\cX$ has compact closure in $\tcX$, so there exists an   accumulation point $y_* \in \overline{\cX} \subset \tcX$ for 
the set $\{ y_n \mid n > 0\} \subset \cX$. We can assume that   $\dT(y_* , y_n) < \delta_*/4$  for all $n > 0$, first by   
passing to a subsequence $\{y_{n_i}\}$ which converges to $y_*$ and satisfies this metric estimate,  and then reindexing the sequence.
 
Let $\ds J_n^x = [y_n - 4\delta_0', y_n+4\delta_0']$, and set   $J_* = [y_* - 3\delta_0', y_* + 3\delta_0']$. 
Then for all $n >0$, we have  $\ds y_n \in (y_* -\delta_0', y_* + \delta_0') \subset J_* \subset J_n^x$.
In particular, $y_1 \in   J_* \subset J_1^x$ is an interior point of $J_*$,     so
$x = \hg_1^x(y_1)$ is an interior point of  $\hg_1^x(J_*)$.

 Also recall from  Proposition~\ref{prop-localhol}, that   $I_n^x = \hg_n^x(J_n^x)$ with 
    $x \in I_n^x$ for all $n$, and   the interval $I_n^x$ has length $|I_n^x| < \delta_0 \exp\{-na/2\} = 8\delta_0' \exp\{-na/2\}$.    
Hence, for $n$ sufficiently large,     the interval $I_n^x$ is contained  in the interior of $\hg_1^x(J_*)$. 
Without loss of generality, we again pass to a 
   subsequence  and reindex the sequence, so that we have  $I_n^x \subset \hg_1^x(J_*)$ and $\ell_{n+1} > \ell_{n}$   for all $n > 0$.
We then    have the inclusions
\begin{equation}\label{eq-inclusions2}
\hg_n^x(J_*) \subset  \hg_n^x(J_n^x) = I_n^x \subset  \hg_1^x(J_*) ~ .
\end{equation}

\begin{figure}[htbp]
\begin{center}
\includegraphics[width=0.7\textwidth]{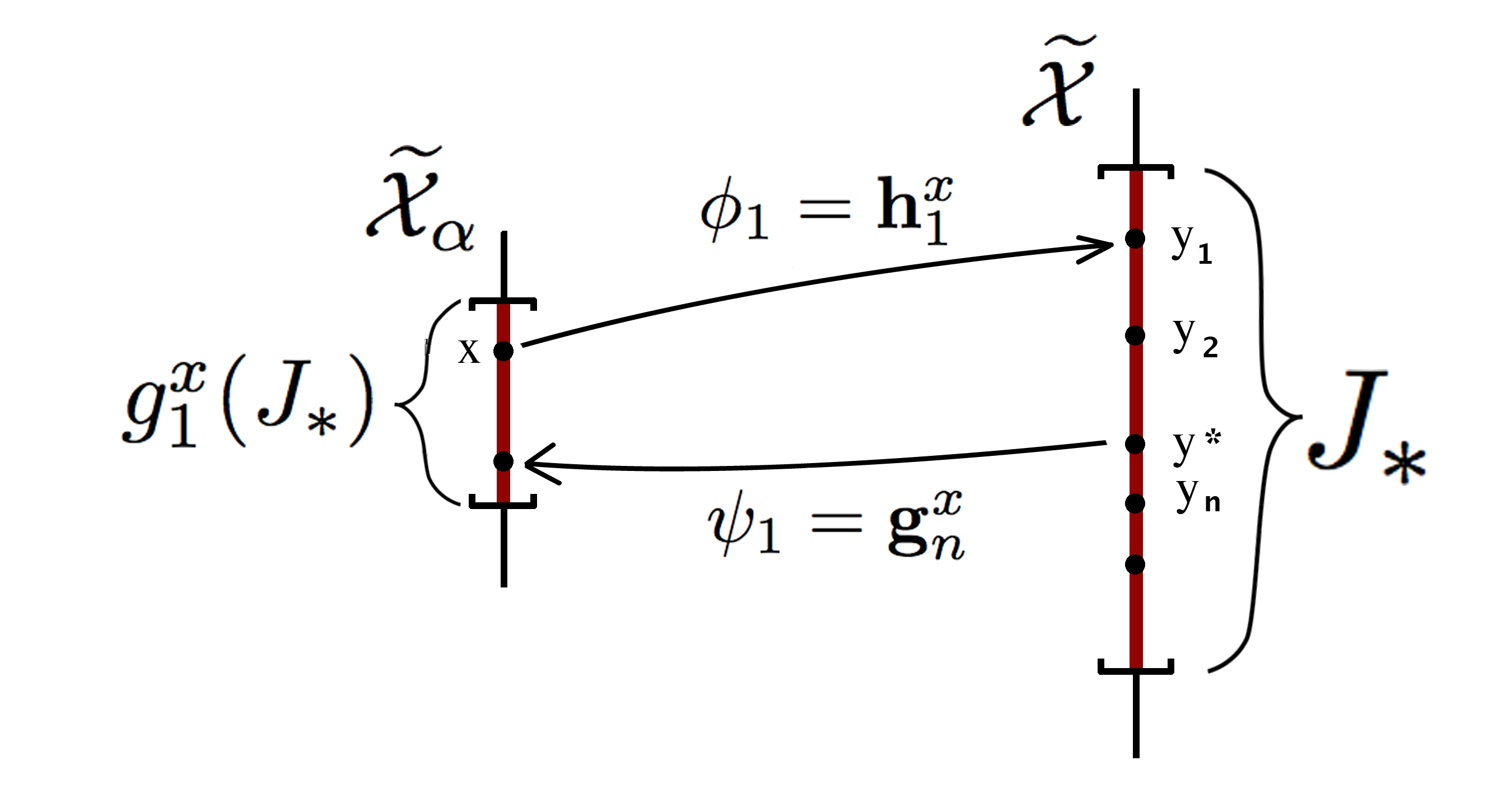}
\caption{The contracting holonomy map $\hh_1^x \circ \hg_n^x$}
\label{Figure5}
\end{center}
\end{figure}

Thus, for each $n > 0$ the composition  $\hh_1^x \circ \hg_n^x   \colon J_* \to   \hh_1^x \circ \hg_1^x(J_*) \subset J_*$ is defined. (See Figure~\ref{Figure5}.)

Recall that $p_1$ denotes the length of the plaque-chain   which defines $\hh_1^x$, and  $C_0$ is the Lipschitz constant defined in the proof of Proposition~\ref{prop-localhol}. Let $N_0$ be chosen so that for  $n \geq N_0$ we have
\begin{eqnarray}
C_0^{p_1} \,  \exp\{-a\, n/2\}  &  < &  \min \, \{\mu , 1/2\}  \label{eq-last1}\\ 
\delta_0' \,  \exp\{-a  \, n/2\}  &  < & \delta_1/2  \label{eq-last2} ~ .
 \end{eqnarray}

 
With the above notations, we then have:
\begin{lemma}\label{lem-hypfp}
Fix $n \geq N_0$, then the map   $\hh_1^x \circ \hg_n^x$ is a hyperbolic contraction on $J_*$ with  fixed-point $v_* \in J_*$ satisfying  $\dT(v_*, y_n) \leq \delta_1/2$ 
and  $(\hh_1^x \circ \hg_n^x)'(v_*) < \mu$. 
\end{lemma}
\proof
By the choice of $C_0$   we have $(\hh_1^x)'(y) \leq C_0^{p_1}$ for all $y$ in its domain.  
Recall that $\hg_n^x$ is the inverse of $\hh_n^x$ which is defined by a plaque-chain of length $\ell_n \geq n$, so the same holds for $\hg_n^x$. 
 The   derivative of  $\hg_n^x$ satisfies the estimates  (\ref{eq-contract1})  by Lemma~\ref{lem-ie-domain},  so we have 
 \begin{equation}\label{eq-contract3}
 \exp\{(-a - 2 \e_1)\, \ell_n\} \leq (\hg_n^x)'(y) \leq    \exp\{(-a + 2\e_1) \, n\} ~. 
\end{equation}
 Thus by \eqref{eq-last1},  for all $y \in J_*$ the composition $\hh_1^x \circ \hg_n^x$  satisfies  
\begin{equation}\label{eq-hyperbolic}
  (\hh_1^x \circ \hg_n^x)'(y) \leq  C_0^{p_1} \,  \exp\{(-a + 2\e_1) \, n\}  <  C_0^{p_1} \,  \exp\{-a\, n/2\}   <  \min \, \{\mu , 1/2\}   
\end{equation}
where we use that the choice of $\e_1 < a/100$ implies that $(-a + 2\e_1) < -a/2$.
Thus,  $\hh_1^x \circ \hg_n^x$ is a hyperbolic contraction on $J_*$ and it follows that $\hh_1^x \circ \hg_n^x$ has a unique fixed-point   $v_* \in J_*$. 
Define a sequence of points  $w_{\ell} = (\hh_1^x \circ \hg_n^x)^{\ell}(y_n) \in J_*$ for $\ell \geq 0$,
then   $\ds v_* = \lim_{\ell \to \infty}~  w_{\ell}$.

 Observe that $\hh_1^x \circ \hg_n^x(y_n) = \hh_1^x(x) =  y_1$, and recall  that  $\dT(y_* , y_n) < \delta_*/4$ for all $n$, hence, $\dT(y_1, y_n) < \delta_*/2$.
Since $w_0 =  y_n$ and $w_1 = y_1$,  
the estimate \eqref{eq-hyperbolic} implies that  
$$\dT(w_{\ell}, w_{\ell +1}) < 2^{-\ell} \cdot \dT(w_0, w_1) < 2^{-\ell} \cdot  \delta_*/2 ~ .$$
 Summing these estimates for $\ell \geq 1$, we obtain   that $\dT(w_0 , v_*)  = \dT(y_n , v_*) \leq \delta_*$  so that 
\begin{equation}\label{eq-estimate7}
\dT(y_*, v_*) \leq  \dT(y_*, y_n) + \dT(y_n, v_*) < 2 \delta_* \leq \delta_1/2 ~ .
\end{equation}
Then  by \eqref{eq-hyperbolic} we have 
$\ds  (\hh_1^x \circ \hg_n^x)'(v_*) \leq \mu$, as was to be shown.
\endproof

The conclusions of Lemma~\ref{lem-hypfp} essentially yield the proof of  Proposition~\ref{prop-fp}, except that it remains to make a change of notation so the results are in the form stated in the proposition, and check that conditions (1) to (5) of Proposition~\ref{prop-fp}.1  are satisfied. This change of notation is done so that the conclusions are in a standard format, which will be invoked recursively in the following Section~\ref{sec-hyperbolic} to prove there exists ``ping-pong'' dynamics in the holonomy pseudogroup $\GF$.

Choose $n \geq N_0$ so that the hypotheses of Lemma~\ref{lem-hypfp} are  satisfied, then define 
  $\phi_1 = \hh_1^x$ and $\psi_1 = \hg_n^x$ so that $\Phi_1 = \phi_1 \circ \psi_1 = \hh_1^x \circ \hg_n^x$, and recall that \quad\quad\quad\quad  $g_1^x(J_*)$
  \begin{equation}\label{eq-inclusions1}
 J_* = [y_* - 3\delta_0', y_* + 3\delta_0'] \subset  J_n^x = [y_n - 4\delta_0', y_n+4\delta_0'] 
\end{equation}
 for $\delta_0'$ and $y_*$ as defined above.     Set $u_1 = v_*$ and $v_1 = \hg_n^x(v_*)$. 
 
 We   check that  conditions  (\ref{prop-fp}.1) and  (\ref{prop-fp}.4) of Proposition~\ref{prop-fp} are satisfied: 
$$\Phi_1(u_1) = \phi_1 \circ \psi_1(u_1) = \hh_1^x \circ \hg_n^x(v_*) = v_* = u_1 ~ , $$
 $$\Psi_1(v_1)  = \psi_1 \circ \phi_1(v_1) = \hg_n^x \circ \hh_1^x(\hg_n^x(v_*)) = \hg_n^x(v_*) = v_1 ~ .$$

Next, for $\cJ_1 = [u_1 - \delta_0', u_1 + \delta_0'] = [v_* - \delta_0', v_* + \delta_0']$ as defined in (\ref{prop-fp}.2),   by the estimate \eqref{eq-estimate7} we have 
$\ds \dT(y_*, v_*)   < 2 \delta_* \leq  \delta_0'/2$
from which it follows that $\cJ_1 \subset J_*$. Then condition (\ref{prop-fp}.3) follows from \eqref{eq-hyperbolic} since $u_1 = v_* \in J_*$. 

 \bigskip
 
Finally, to show condition (\ref{prop-fp}.5) of Proposition~\ref{prop-fp} is satisfied,   recall that $\psi_1(y_n) = \hg_1^x(y_n) =  x$,        that  $\dT(y_n , v_*) <  \delta_* \leq 1$  by the proof of  Lemma~\ref{lem-hypfp}, and that  $\delta_* = \min\{1,\delta_0'/4,\delta_1/4\}$. 
  Also, the estimate    (\ref{eq-contract1})   combined with  \eqref{eq-last2}  and the choice of  $\delta_0' \leq 1$ in Definition~\ref{def-delta0} yields that, for all $y \in J_*$
\begin{equation}\label{eq-contraction3}
 (\hg_n^x)'(y) \leq    \exp\{(-a + 2\e_1) \, \ell_n\}  <  \delta_0' \cdot  \exp\{-a\, n/2\}   <  \delta_1/2 ~ .
\end{equation}
Thus, by the   Mean Value Theorem and the estimate $\dT(y_n, v_*)  \leq \delta_* \leq 1$, we have that  
$$\dT(x, v_1) = \dT(\hg_n^x(y_n), \hg_n^x(v_*)) \leq \delta_1/2 \cdot \dT(y_n, v_*) \leq \delta_1/2 ~ .$$
For any $y \in \cJ_1 = [v_* - \delta_0', v_* + \delta_0']$ we also have that
   $$\dT(\hg_n^x(y), v_1) = \dT(\hg_n^x(y), \hg_n^x(v_*)) \leq \delta_1/2 \cdot \dT(y, v_*) \leq \delta_0' \delta_1/2 < \delta_1/2 ~.$$
Thus,   
$$ \dT(\hg_n^x(y), x) \leq \dT(\hg_n^x(y), v_1) + \dT(x, v_1) <  \delta_1 ~ , $$ 
so that  $\cK_1 =\psi_1(\cJ_1)  \subset [x - \delta_1, x+ \delta_1]$, as was to shown.
 
   This   completes the proof of Proposition~\ref{prop-fp}. 
\endproof

\section{Hyperbolic sets with positive measure}\label{sec-hyperbolic}

 The main result of this section is:
 \begin{thm}\label{thm-entropy} Let $\F$ be a  $C^1$-foliation of codimension-one of a compact manifold $M$ for which 
 $\E$ has positive Lebesgue measure. Then $\F$ has a hyperbolic resilient leaf, and hence  the geometric entropy  $\h > 0$. 
\end{thm}
 
The assumption that the Lebesgue measure $|\E| > 0$ is used in two ways. First, 
the set $\E$ is an increasing union  of the sets $\Ea$   for   $a > 0$, so $|\E| > 0$  implies   $|\Ea| > 0$      for some $a> 0$. For each $x \in \Ea$,  we obtain from Proposition~\ref{prop-fp} uniform hyperbolic contractions with  fixed-points arbitrarily close to the given $x \in E$, and with prescribed bounds on their domains.
 
  Secondly, almost every point of a measurable set is a point of positive Lebesgue density, hence  $|\Ea| > 0$ implies that 
$\Ea$ has a ``pre-perfect'' subset of points with expansion greater than $a$. This observation enables us to   construct an infinite sequence of hyperbolic fixed-points arbitrarily close to the support of $\Ea$, whose  domains have to eventually overlap since the closure $\overline{\cX}$ is compact. This yields the existence of  a resilient orbit for $\GF$, hence a   ping-pong game dynamics as defined in Section~\ref{subsec-resilient},  which implies that $\h > 0$.

\begin{defn}
A set $\cE$ is said to be \emph{pre-perfect} if it is non-empty, and its closure $\overline{\cE}$ is a perfect set. Equivalently, $\cE$ is pre-perfect if it is not empty, and no point is isolated.  
\end{defn}
The following observation is a standard property of sets with positive Lebesgue measure.
 \begin{lemma}\label{lem-preperfect}
If $X \subset \mR^q$ has positive Lebesgue measure, then there is a pre-perfect subset $\cE \subset X$.
 \end{lemma}
\proof
 Let $\cE \subset X$ be the set of points with  Lebesgue density $1$. Recall that this means that  for each $x \in X$ and each $\delta > 0$, 
 the Lebesgue measure  $|B_X(x,\delta) \cap X| > 0$,  and   $\ds \lim_{\delta \to 0}\frac{ |B_X(x,\delta) \cap X|}{|B_X(x,\delta)|} = 1$. 
 
 It is a standard fact of Lebesgue measure theory that $|\cE| = |X|$, so that 
 $|X| > 0$ implies that $\cE \ne \emptyset$. Moreover, if $x \in \cE$ is isolated in $\cE$, then $x$ is a point with Lebesgue density $0$, thus each  $x \in \cE$ cannot be isolated. It follows that $\cE$ is pre-perfect.
\endproof

Theorem~\ref{thm-entropy} now follows from Lemma~\ref{lem-preperfect} and the following result:

\begin{prop}\label{prop-resilient}
Let $a>0$, and suppose  there exists a pre-perfect subset  $\cE \subset \Ea$, 
then $\F$ has a resilient leaf contained in the closure $\overline{\Ea}$.
\end{prop}
\proof
 Let $a > 0$ and let $\cE \subset \Ea$ be a pre-perfect set.  The saturation of a pre-perfect set under the action of the holonomy pseudogroup $\GF$ is pre-perfect, so we can assume that $\cE$ is saturated.
   We assume that $\F$ does not have a resilient leaf in   $\overline{\Ea}$, and show this leads to  a contradiction. 

   We follow the notation introduced in the proof of Proposition~\ref{prop-fp}, which will be invoked repeatedly, and the resulting maps and constants will be labeled  according to the stage of the induction.
Choose   $0 < \e_1 < \min\{\e_0, a/100\}$, and let $\delta_0$ be chosen as in Definition~\ref{def-delta0}.

   Fix a choice of $0 < \mu < 1$, and choose  $0 < \delta_1 < \delta_0$ and  $x_1 \in \cE \cap \cX_{\alpha}$.  
Then by Proposition~\ref{prop-fp},     there exists holonomy maps $\phi_1, \psi_1 \in \GF$ and points $u_1  \in \overline{\cX}$ and $v_1 = \psi_1(u_1)$, such that $\dT(x_1,v_1) < \delta_1$ and which are fixed-points for the maps  $\Phi_1$, $\Psi_1$ respectively.   Moreover, we have the sets
 \begin{enumerate}
\item  $\cJ_1 \equiv [u_1 - \delta_0, u_1+ \delta_0]$ 
\item  $\cI_1 \equiv \Phi_1(\cJ_1) \subset (u_1 - \delta_0, u_1+ \delta_0)$
\item  $\cK_1 \equiv \psi_1(\cJ_1) \subset (x_1 - \delta_1, x_1+ \delta_1)$
\end{enumerate}
 whose properties were given in Proposition~\ref{prop-fp}. In particular, $\Phi_1 \colon \cJ_1 \to \cI_1 \subset \cJ_1$ is a hyperbolic contraction with fixed-point $u_1$. In particular, note that $\ds \bigcap_{\ell > 0} \, \Phi_1^{\ell}(\cJ_1) = \{ u_1 \}$.

If the orbit of $u_1$ under $\GF$ intersects $\cJ_1$   in   a point other than $u_1$, 
then by definition, $u_1$ is a hyperbolic resilient point, which by assumption does not exist. Therefore, the $\GF$-orbit of $u_1$ intersects the interval $\cJ_1$ exactly in  the interior point $u_1$, and    intersects     $\cK_1$ exactly in  the interior point $v_1$. 

 Note that $x_1 \in \cK_1 \cap \cE$ so there exists $x_2 \in (\cK_1 - \{x_1, v_1\}) \cap \cE$ as $\cE$ is pre-perfect. Choose $0 < \delta_2 < \delta_1$ so that 
 $\ds (x_2 -  \delta_2, x_2 + \delta_2) \subset (\cK_1 - \{x_1, v_1\})$. The $\GF$-orbit of $v_1$   intersects $\cK_1$ only in the point $v_1$, thus the interval $(x_2 -  \delta_2, x_2 + \delta_2)$ is disjoint from the orbit of  $v_1$. We then repeat the construction in the proof of Proposition~\ref{prop-fp},    to obtain  holonomy maps $\phi_2, \psi_2 \in \GF$ and points $u_2  \in \overline{\cX}$ and $v_2 = \psi_2(u_2)$, such that $\dT(x_2,v_2) < \delta_2$ and which are fixed-points for the maps  $\Phi_2$, $\Psi_2$ respectively. Again, define the sets
 \begin{enumerate}
\item  $\cJ_2 \equiv [u_2 - \delta_0, u_2+ \delta_0]$ 
\item  $\cI_2 \equiv \Phi_2(\cJ_2) \subset (u_2 - \delta_0, u_2+ \delta_0)$
\item  $\cK_2 \equiv \psi_2(\cJ_2) \subset [x_2 - \delta_2, x_2+ \delta_2]$ .
\end{enumerate}

We then repeat this construction recursively. Let $\{u_1, u_2, \ldots\} \subset \overline{\cX}$ be the resulting centers of contraction for the hyperbolic maps $\{\Phi_{i} \mid i > 0\}$. As $\overline{\cX}$ is compact, there exists an accumulation point $u_* \in \overline{\cX}$. In particular, there exists distinct indices $i_1, i_2 > 0$ such that $\dT(u_*, u_{i_1}) < \delta_0/10$ and $\dT(u_*, u_{i_2}) < \delta_0/10$ and thus $\dT(u_{i_1}, u_{i_2}) < \delta_0/5$. 

 Recall that the intervals $\cJ_{i_1} = [u_{i_1} - \delta_0, u_{i_1}+ \delta_0]$ and $\cJ_{i_2} = [u_{i_2} - \delta_0, u_{i_2}+ \delta_0]$ have uniform width, and moreover  
 $\ds \{u_{i_1}, u_{i_2}\} \subset \cJ_{i_1} \cap \cJ_{i_2}$. 
 As $u_{i_1}$ and $u_{i_2}$ are disjoint fixed-points of hyperbolic attractors, we can choose   integers $m_1, m_2 > 0$ so that 
 $\ds \Phi_{i_1}^{m_1}(\cJ_{i_1}) \cap \Phi_{i_2}^{m_2}(\cJ_{i_2}) = \emptyset$ and $\Phi_{i_j}^{m_j}(\cJ_{i_j}) \subset \cJ = \cJ_{i_1} \cap \cJ_{i_2}$ for $j=1,2$. 
Then  the action of the contracting maps  $\bH = \Phi_{i_1}^{m_1}$ and $\bG = \Phi_{i_2}^{m_2}$ on $\cJ$ define  a ``ping-pong game'' as in Definition~\ref{defn-ppg}.

  Now let $x = u_{i_1}$,  $y = \bG(x) \ne x$, then $\bH^{\ell}(y) \to x$ as $\ell \to \infty$, so that the orbit of $x$ under the action $\GF$ is resilient,  contrary to assumption. 
 
     Hence, if there exists a pre-perfect set $\cE \subset \Ea$ for  $a > 0$,  then there   exists a resilient leaf.
\endproof

\section{Open manifolds}\label{sec-open}

 In this section, we   extend the methods above from compact manifolds to open  manifolds, using the techniques of \cite[Section~5]{Hurder1986}. 
 \begin{thm}\label{thm-openGV} 
 Let $\F$ be a codimension-one  $C^2$-foliation of an open complete manifold $M$. If the   Godbillon-Vey class  $GV(\F) \in H^3(M; \mR)$ is non-zero, the $\F$  has  a hyperbolic resilient leaf. 
 \end{thm}
\proof
 The class $GV(\F) \in H^3(M; \mR)$  is   determined by its pairing with the compactly supported cohomology group $H^{m-3}_c(M; \mR)$, so $GV(\F) \ne 0$ implies there exists a closed $m-3$ form $\xi$ with compact support on $M$ such that $\langle GV(\F), [\xi]\rangle \not= 0$. Let $|\xi| \subset M$ denote the  support of $\xi$, which is a compact set. 
As the support $|\xi|$ is compact, there is a  finite open cover of $|\xi|$ by a regular foliation atlas $\{ (U_{\alpha},\phi_{\alpha})  \mid    \alpha \in \cA \}$ for $\F$ on $M$ (as in Section~\ref{sec-basics} above).  Let $M_0$ denote the union of the sets 
$\{ U_{\alpha}   \mid    \alpha \in \cA \}$, then   the closure $\overline{M_0}$ is a compact subset of $M$ and $|\xi| \subset M_0$. Thus we have   $GV(\F|M_0) \not= 0$. If $M_0$ is not connected, we can choose a connected component $M_1 \subset M_0$ for which $GV(\F|M_1) \not= 0$. Thus, we may assume that $M_0$ is connected.

The proof of Theorem~\ref{thm3} used only the properties of the pseudogroup generated by a regular foliation atlas $\{ (U_{\alpha},\phi_{\alpha})  \mid    \alpha \in \cA \}$ --  the compactness of $M$ was not used except in the construction of this atlas. The definition and properties of the Godbillon measure also apply to open manifolds, as was discussed in  \cite[Section 5]{Hurder1986}. Hence, by the same proof  we obtain that the set ${\rm E}(\F|M_0)$ has positive measure. 

The proofs of Propositions~\ref{prop-fp} and \ref{prop-resilient} use only   the assumption that    the pseudogroup $\GF$ is compactly generated, as defined by Haefliger \cite{Haefliger2002}, and   do not require the compactness of $M$, hence apply directly to show that ${\cG}\F|M_0$ has a hyperbolic resilient point if ${\rm E}(\F|M_0)$ has 
positive measure. Thus, $\F|M_0$ must have a resilient leaf, and so also must $\F$.
\endproof

 Here is an application of Theorem~\ref{thm-openGV}. 
 Let $\bB\G_1^{(2)}$ denote the Haefliger classifying space of codimension--one $C^2$-foliations \cite{Haefliger1970,Haefliger1971}. There is a universal Godbillon-Vey class $GV \in H^3(\bB\G_1^{(2)}; \mR)$ 
such that for every codimension--one $C^2$-foliation $\F$ of a manifold $M$, there is a classifying map $h_{\F} \colon M \to \bB\G_1^{(2)}$ such that $h_{\F}^* GV = GV(\F)$   (see \cite{BottHaefliger1972,Lawson1977}.)  The first two integral homotopy groups 
$\pi_1(\bB\G_1^{(2)})  = 0 =  \pi_2(\bB\G_1^{(2)})$, while Thurston showed in \cite{Thurston1972} that the Godbillon-Vey class defines a \emph{surjection} $GV \colon \pi_3(\bB\G_1^{(2)}) \to \mR$.  It follows from Thurston's work in \cite{Thurston1976}, that for a closed oriented $3$-manifold $M$ and any $a >0$,  there exists a  codimension--one  foliation $\F_a$ on $M$  such that  $\langle GV(\F_a) , [M]\rangle = a$. Each such foliation $\F_a$ for $a \ne 0$ must then have resilient leaves.

More generally, given any finite CW complex $X$, a continuous map $h \colon X \to \bB\G_1^{(2)}$ defines a foliated microbundle over $X$, whose total space $M$ is an open manifold with a codimension--one foliation $\F_h$ such that $h^*GV = GV(\F_h)$.  This is discussed in detail by Haefliger \cite{Haefliger1971}, who introduced the technique. (See also Lawson \cite{Lawson1977}.)  Thus, using homotopy methods to construct the map $h$ so that $h^*GV \not= 0$, one can construct many examples  of open foliated manifolds with non-trivial Godbillon-Vey classes.
Theorem~\ref{thm-openGV} implies that all such examples have  resilient leaves.


\small


 \bigskip

\end{document}